\documentclass[11pt]{article}    
\usepackage{epsfig,amsmath,amssymb,multirow,color,tabularx}     

\begin{document}

\title{A matched alternating direction implicit (ADI) method for solving
the heat equation with interfaces}

\author{Shan Zhao$^{1,2,}$\footnote{
Corresponding author. Tel: 1-205-3485155, Fax: 1-205-3487067,
Email: szhao@ua.edu}\\
$^1$ Department of Mathematics,
University of Alabama, Tuscaloosa, AL 35487, USA. \\
$^2$ Beijing Computational Science Research Center, 
Beijing 100084, PR China.
}

\date{\today} 
\maketitle

\begin{abstract}
A novel Douglas alternating direction implicit (ADI) method is proposed in this
work to solve a two-dimensional (2D) heat equation with interfaces. 
The ADI scheme is a powerful finite difference method for solving parabolic equations,
due to its unconditional stability and high efficiency. 
However, it suffers from a serious accuracy reduction in space  for interface problems
with different materials and nonsmooth solutions.
If the jumps in a function and its derivatives are known across the interface, 
rigorous ADI schemes have been successfully constructed in the literature 
based on the immersed interface 
method (IIM) so that the spatial accuracy can be restored. 
Nevertheless, the development of accurate and stable ADI methods 
for general parabolic interface problems with 
physical interface conditions that describe jumps of a function and its flux,
remains unsolved. 
To overcome this difficulty, a novel tensor product decomposition is proposed
in this paper to decouple 2D jump conditions into essentially one-dimensional
(1D) ones. These 1D conditions can then be incorporated into the ADI 
central difference discretization, using the matched interface and boundary
(MIB) technique. 
Fast algebraic solvers for perturbed tridiagonal
systems are developed to maintain the computational efficiency. 
Stability analysis is conducted through eigenvalue spectrum analysis, which
numerically demonstrates the unconditional stability of the proposed ADI method. 
The matched ADI scheme achieves the first order of accuracy in time 
and second order of accuracy in space 
in all tested parabolic interface problems with
complex geometries and spatial-temporal dependent jump conditions. 

\noindent {\bf Keyword:} 
Heat equation; Parabolic interface problem;
Jump conditions;
Alternating direction implicit (ADI);
Matched interface and boundary (MIB).

\noindent {\bf MSC:}
65M06, 
65M12, 
35K05. 
\end{abstract}

\section{Introduction}

In this paper, we propose a new alternating direction implicit (ADI) method for solving
two-dimensional (2D) parabolic interface problems
\begin{equation}\label{heat}
\frac{\partial u}{\partial t} = \nabla \cdot ( \alpha \nabla u) + f, \quad \mbox{in}~~ 
\Omega=\Omega^- \cup \Omega^+
\end{equation}
where $u(x,y,t)$ is a function of interests, e.g. the temperature, 
$\alpha$ is the diffusion coefficient and $f(x,y,t)$ is a source.
For simplicity, the  domain 
$\Omega$ is assumed to be a rectangle one, with
proper boundary conditions prescribed for $u$ on  $\partial \Omega$. 
Across a material interface $\Gamma$ separating two media $\Omega^-$ and
$\Omega^+$, the diffusion coefficient $\alpha$ is discontinuous, while the source term
$f(x,y,t)$ may be even singular. 
Physically,  the solution $u$ on both sides of $\Gamma$ 
is related analytically via jump conditions
\begin{equation}\label{jump}
[u] = u^+ - u^- = \phi(s,t),  \quad
[\alpha u_n] = \alpha^+ \frac{\partial u^+}{\partial n} 
- \alpha^- \frac{\partial u^-}{\partial n}  = \psi(s,t),
\end{equation}
where  $s$ is the arc-length parametrization
of the interface $\Gamma$, and $n$ the unit normal direction. 
The superscript, $-$ or $+$, denotes the limiting value of a function from one side or
the other of the interface. 
We note that (\ref{jump}) takes a quite general form, while for many applications, 
we have simply $\phi=\psi=0$. 

The parabolic interface problem governed by (\ref{heat}) and (\ref{jump})
appears in many physical and engineering applications, 
such as the continuous casting in the metallurgical industry, 
the freezing process of perishable foodstuffs in the food engineering, 
and the magnetic fluid hyperthermia treatment of cancer. 
The analysis of conductive heat transfer process over composite media is
indispensable in these applications. 

Since the physical solution is non-smooth
or even discontinuous across the interface, the standard numerical methods
often perform poorly for the parabolic interface problem. 
To restore the accuracy near the interface, the jump conditions (\ref{jump}) 
have to be incorporated into the numerical discretization in certain manner. 
For finite element methods 
\cite{ChenZou,Sinha05,Sinha09,Attanayake}
and finite volume methods \cite{Wang10},
many rigorous interface treatments have  been proposed to deliver high  accuracy
in solving parabolic interface problems.

As one of the most successful finite difference methods for solving material
interface problems, the immersed interface method (IIM) was originally introduced 
by LeVeque and Li \cite{LevLi} for solving elliptic equations with discontinuous
coefficients and singular sources. 
By rigorously imposing jump conditions via local Taylor expansions,
the IIM achieves the second order of accuracy for complex elliptic interfaces. 
The development of the IIM for solving parabolic equations has been considered in 
\cite{Kandilarov04,Kandilarov07,Bouchon10a,Bouchon10b,AdamsLi}.
In \cite{Kandilarov04,Kandilarov07}, the authors constructed the IIM together with the
implicit Euler time integration for parabolic equations with singular own sources.
The IIM has been applied in \cite{Bouchon10a,Bouchon10b} for solving the Poisson 
equation over moving irregular domains or 
over fixed domains with the Neumann boundary condition. 
The maximum principle preserving IIM has been proposed in \cite{AdamsLi} by
Adams and Li for solving the convection-diffusion equation with general jump conditions.
In this scheme,  the advection term is discretized explicitly while the 
diffusion term is treated implicitly. A fast multigrid method is implemented to 
efficiently solve the linear system of equations for the implicit time stepping  \cite{AdamsLi}. 

As one of the most successful finite difference methods for solving parabolic
equations, the classical ADI method \cite{DougPeace,Douglas,PeaceRach} can
be written as some perturbations of multidimensional implicit methods, such as the 
Crank-Nicolson and backward Euler. In general, the ADI method is unconditionally stable 
for parabolic problems without interfaces so that a large time increment is admissible, 
which in turn can lead to a faster simulation.
The major attraction of the ADI method, as compared with other implicit methods,
is that it reduces a multidimensional problem
to sets of independent one-dimensional (1D) problems of tridiagonal structures,
and such matrices can be efficiently solved using the  Thomas algorithm \cite{ADIbook}.
Moreover, the essentially 1D feature of the ADI computations allows 
a tremendously efficient parallelization, 
including on modern Graphics Processing Units (GPUs) \cite{Tay,Wei13}.
Neglecting jump conditions, the classical ADI method has been applied to 
solve interface problems with sharp or smeared interfaces, see for example 
\cite{Bates09,Chen11,Geng13,Zhao14,Tian14}. 
However, for sharp interface problems, the central finite difference approximation
degrades to the first order of accuracy in space
and the precious unconditional stability may be lost \cite{Geng13,Zhao14}. 
Therefore, it is highly desired to construct novel ADI methods 
for solving parabolic interface problems
without compromising stability and spatial accuracy. 

The first rigorous interface treatment in the ADI framework was due to
Li and Mayo  \cite{LiMayo}, in which a homogeneous 2D heat equation 
with a constant diffusion coefficient and a singular source is solved. 
For such an interface setting, 
the jump conditions (\ref{jump}) take a simpler form with $\alpha=1$. Thus,
higher order jump conditions can be simply derived \cite{LiMayo}. 
This allows the construction of a second order accurate IIM-ADI method by 
adding some correction terms into the classical ADI scheme
for irregular points near the interface.
Considering the same jump conditions, Liu and Zheng have extended
the IIM-ADI method to solve a 2D homogeneous 
convection-diffusion equation \cite{Liu13}
and a three-dimensional (3D) homogeneous heat equation \cite{Liu14}. 

However, grand difficulties are encountered when the IIM-ADI method
\cite{LiMayo,Liu13,Liu14}  is generalized in \cite{LiShen}
to solve a 2D heat equation with nonhomogeneous media,
i.e., $\alpha$ being a piecewise constant. 
By still assuming simple jump conditions with prescribed function and
derivative jumps, i.e., given $[u]$ and $[u_n]$ values, 
the second order jump condition for $[u_{nn}]$ cannot be simply derived
from the governing equation now. 
A rather complicated second order jump condition is considered in  \cite{LiShen}
so that the SOR iterative method has to be used to solve a 2D linear 
system in each step of the Crank-Nicolson time integration.
We note that, on the other hand, if the general jump conditions like (\ref{jump})
are used, the second order jump conditions could be naturally derived. 
Nevertheless, the construction of the IIM-ADI scheme through introducing 
correction terms remains a challenge for general jump conditions. 
Therefore, the development of accurate and stable ADI methods for 
the  parabolic interface problem  (\ref{heat}) and (\ref{jump})
is essentially an open problem. 

The objective of this paper is to propose a novel matched ADI method to
overcome the aforementioned difficulties for solving 
general parabolic interface problems. 
The proposed matched ADI method is formulated based on our previous 
interface scheme,  the matched interface and boundary (MIB) method,
originally developed for solving elliptic and hyperbolic interface problems
\cite{Zhao04,Zhou06}. 
One distinction between the MIB and IIM is that the MIB interface modeling 
just needs zeroth and first order jump conditions, i.e., (\ref{jump}) so that the difficulty
associated with the second order jump conditions of the IIM is simply 
bypassed in the MIB approach. 
However, the MIB scheme has never been applied to a parabolic interface
problem before.
Moreover, the existing MIB scheme cannot be directly utilized in the ADI
formulation, because the 2D MIB interface treatment will couple $x$ and $y$
directions simultaneously. 
The most significant contribution of this work is the introduction
of a novel tensor-product
decomposition of jump conditions  (\ref{jump}), which decouples  2D
jump conditions into 1D ones, in the same spirit of the ADI method. 
Then, 1D MIB interface treatments will be developed  in space to secure 
a second order of accuracy. 
Fast algebraic solvers based on the Thomas algorithm will be developed
to solve 1D linear systems efficiently. 
The stability proof of the matched ADI algorithm is highly non-trivial,
because the finite difference weights of the MIB discretization 
depend on the interface geometry in an unpredictable manner.
In the present study, through calculating the spectral radius,
the proposed matched ADI method is numerically verified to be unconditionally stable.

The rest of this paper is organized as follows. 
Section \ref{Sec:theory} is devoted to the theory and algorithm
of the proposed matched ADI method. 
Numerical tests are carried out to validate the proposed method by 
considering various particular forms for the jump conditions (\ref{jump}).
Finally, a conclusion ends this paper. 

\section{Theory and Algorithm}\label{Sec:theory}

Consider an interface problem, in which  $\Omega^-$ is interior to $\Omega^+$.
Define a uniform mesh partition of the computational domain $\Omega$.
Without the loss of generality, we assume that the grid spacing $h$ in
both $x$ and $y$ directions is the same and 
one grid line cuts the interface $\Gamma$ at most twice.
Denote the time increment to be $\Delta t$ and take
$N_x$ and $N_y$ as the number of grid points in each direction.
To facilitate the following discussions, we adopt a
notation at node $(x_i,y_j,t_k)$:
$u^k_{i,j}=u(x_i,y_j,t_k)$.

\subsection{Temporal discretization}
We first rewrite the heat equation (\ref{heat2}) by dividing $\alpha$ throughout
\begin{equation}\label{heat2}
\frac{1}{\alpha}  \frac{\partial u}{ \partial t}
=\frac{\partial^2 u }{\partial x^2} + \frac{\partial^2 u }{\partial y^2}
+\frac{f}{\alpha},   \quad \mbox{in}~~ \Omega^- ~~\mbox{or}~~ \Omega^+,
\end{equation}
If the jump conditions (\ref{jump}) are rigorously enforced in the numerical
discretization, the numerical solution to (\ref{heat2}) will be identical to 
that of (\ref{heat}), whereas 
(\ref{heat2}) allows an easier formulation for the ADI method.
The semi-discretization of  (\ref{heat2}) using the implicit Euler time
integration at a general spacial node $(x_i,y_j)$ reads
\begin{equation}\label{ADIsemi}
\frac{u^{k+1}_{i,j} - u^k_{i,j} }{\alpha \Delta t}
=\delta_{xx} u^{k+1}_{i,j} + \delta_{yy} u^{k+1}_{i,j} 
+\frac{f^{k+1}_{i,j}}{\alpha},
\end{equation}
which is first order accurate in time.
Here $\delta_{xx}$ and $\delta_{yy}$ are discrete operators for
finite difference approximations in $x$ and $y$ directions. 
We propose a first order Douglas ADI method for solving (\ref{heat2}),
\begin{align}
(\frac{1}{\alpha} -  \Delta t \delta_{xx}) u^*_{i,j} & = 
(\frac{1}{\alpha} + \Delta t \delta_{yy}) u^k_{i,j}
+ \frac{\Delta t}{\alpha} f^{k+1}_{i,j}, \nonumber \\
(\frac{1}{\alpha} -  \Delta t \delta_{yy}) u^{k+1}_{i,j} & = 
\frac{1}{\alpha} u^*_{i,j} - \Delta t \delta_{yy} u^k_{i,j}. \label{ADI} 
\end{align}
To see the connection between (\ref{ADIsemi}) and (\ref{ADI}), we can 
eliminate $u^*_{i,j}$ in (\ref{ADI}),
\begin{equation}\label{Doug2}
(\frac{1}{\alpha}-\Delta t \delta_{xx}) 
(\frac{1}{\alpha}-\Delta t \delta_{yy}) u^{k+1}_{i,j}
= (\frac{1}{\alpha} +\Delta t \delta_{yy}) \frac{1}{\alpha} u^k_{i,j}
- (\frac{1}{\alpha}-\Delta t \delta_{xx}) \Delta t \delta_{yy} u^k_{i,j}
+ \frac{\Delta t}{\alpha^2} f^{k+1}_{i,j}.
\end{equation}
After fully expanding terms,  Eq. (\ref{Doug2}) can be written into the form
\begin{equation}\label{Doug3}
(\frac{1}{\alpha}-\Delta t \delta_{xx} -\Delta t \delta_{yy}
+ \alpha \Delta t^2 \delta_{xx} \delta_{yy} ) u^{k+1}_{i,j}
= (\frac{1}{\alpha} + \alpha \Delta t^2 \delta_{xx} \delta_{yy})u^{k}_{i,j}
+ \frac{\Delta t}{\alpha} f^{k+1}_{i,j}.
\end{equation}
If we drop the higher order perturbation term 
$\alpha \Delta t^2 \delta_{xx} \delta_{yy} u$ on both hand sides of
(\ref{Doug3}), we actually obtain an equivalent form of the implicit
Euler scheme (\ref{ADIsemi})
\begin{equation}\label{imEuler}
(\frac{1}{\alpha}-\Delta t \delta_{xx} -\Delta t \delta_{yy}) u^{k+1}_{i,j}
= \frac{1}{\alpha} u^{k}_{i,j} + \frac{\Delta t}{\alpha} f^{k+1}_{i,j}.
\end{equation}
In this work, the Douglas scheme (\ref{ADI}) will be employed in all 
ADI computations, while the implicit Euler scheme  (\ref{imEuler}) can be 
used in the theoretical analysis.

\subsection{Spatial discretization}
We next consider the spatial discretization. 
For nodes away from the interface, a central difference approximation
is used, e.g.,
\begin{equation}\label{uyy}
\delta_{yy} u^{k+1}_{i,j}
:=\frac{1}{h^2} (u^{k+1}_{i,j-1} - 2 u^{k+1}_{i,j} + u^{k+1}_{i,j+1}).
\end{equation}
In the present study, a proper Dirichlet or Neumann 
boundary condition for $u$ is assumed to be given 
on the boundary $\partial \Omega$. Such a boundary condition is 
implemented as in the classical ADI schemes 
\cite{DougPeace,Douglas,PeaceRach}.
If a nontrivial boundary condition is encountered, the advanced 
MIB boundary closure method \cite{Zhao07,Zhao09}
can be utilized to enclose such a 
condition into the finite difference discretization.

For nodes near the interface $\Gamma$, novel interface treatments will
be developed to correct discrete finite difference operators 
$\delta_{xx}$ and $\delta_{yy}$  via rigorously imposing the jump
conditions (\ref{jump}). 
To this end, we consider some tensor product decompositions of jump
conditions (\ref{jump}) in the ADI framework. 
At an interface point, we denote the outer normal and tangential directions
as $n$ and $\tau$, respectively. 
Denote the angle between $n$ and the $x$-axis as $\theta$. 
Coordinate transformations can be employed to convert between the
derivatives 
\begin{equation}
\frac{\partial}{\partial n} = \cos \theta \frac{\partial}{\partial x} 
+ \sin \theta \frac{\partial}{\partial y}, \quad
 \frac{\partial}{\partial \tau} = -\sin \theta \frac{\partial}{\partial x} 
+ \cos \theta \frac{\partial}{\partial y}.
\end{equation}
Based on (\ref{jump}), one more jump condition can be derived by
differentiating along the interface $\Gamma$:
$[u_{\tau}] =  \frac{\partial \phi}{\partial \tau} = \phi_{\tau}$.
We thus have three zeroth and first order jump conditions
\begin{equation}\label{jump2}
[u] = \phi, \quad [u_{\tau}] = \phi_{\tau}, \quad
[\alpha u_n] = \psi. 
\end{equation}
However, jump conditions given in (\ref{jump2}) cannot be applied in a 1D
manner in the ADI algorithm. To illustrate this, we transform the flux
jump condition into Cartesian directions
\begin{equation}\label{jump-xy}
\cos \theta [\alpha u_x] + \sin \theta [\alpha u_y] = \psi.
\end{equation}
A strong coupling in $x$ and $y$ directions is clearly seen. 

In this paper, we propose to decompose the 2D jump conditions (\ref{jump2})
into two sets of essentially 1D jump conditions. 
We illustrate the idea by considering the $x$ direction formulation. 
The $y$ direction can be similarly treated. 
Consider an interface point $(x_{\Gamma},y_j)$ which is the
intersection point between one $x$ grid line and $\Gamma$. If the normal
direction $n$ at this point happens to be along the $x$ direction, we have simply
1D jump conditions 
$[u]=\phi$ and $[\alpha u_x]= \pm \psi$. For more general scenario, in which
$n$ is not along the $x$ direction, we will analytically derive a hybrid jump condition
using $x$ and $\tau$. 
In particular, we have $u_y=\sec \theta u_{\tau} + \tan \theta u_x$. By
substituting this  into (\ref{jump-xy}), we arrive at 
\begin{equation}\label{jump-xtau}
[\alpha u_x] + \sin \theta [\alpha u_{\tau}] = \cos \theta \psi.
\end{equation}
Even though Eq. (\ref{jump-xtau}) looks similar to  Eq. (\ref{jump-xy}), 
Eq. (\ref{jump-xtau}) allows a 1D numerical approximation, whereas 
Eq. (\ref{jump-xy}) dose not. This is because the jump value of $[u_{\tau}]$ is known,
while that of $[u_y]$ is not.
In fact, Eq. (\ref{jump-xtau}) can be further rewritten as
\begin{equation}\label{jump-x}
[\alpha u_x] = \cos \theta \psi - \sin \theta (\alpha^+ - \alpha^-) u^+_{\tau}
- \sin \theta \alpha^- \phi_{\tau} = \bar{\psi},
\end{equation}
where $u^+_{\tau}$ shall be evaluated based on finite difference approximations 
using grid nodes exclusively from the positive side of $\Gamma$, i.e., $\Omega^+$. 
Once  $u^+_{\tau}$ is accurately estimated, $\bar{\psi}$ is known. 
We thus have essentially 1D jump conditions
$[u]=\phi$ and $[\alpha u_x]=\bar{\psi}$. Similarly, we derive the following
essentially 1D jump conditions in the $y$ direction,
\begin{equation}\label{jump-y}
[u]=\phi, \quad 
[\alpha u_y] = \sin \theta \psi + \cos \theta (\alpha^+ - \alpha^-) u^+_{\tau}
+ \cos \theta \alpha^- \phi_{\tau} = \hat{\psi}.
\end{equation}
We note that $\bar{\psi}$ and $\hat{\psi}$ can also be evaluated through calculating
$u^-_{\tau}$ from the negative side of the interface.

\begin{figure*}[!tb]
\begin{center}
\begin{tabular}{ccc}
\psfig{figure=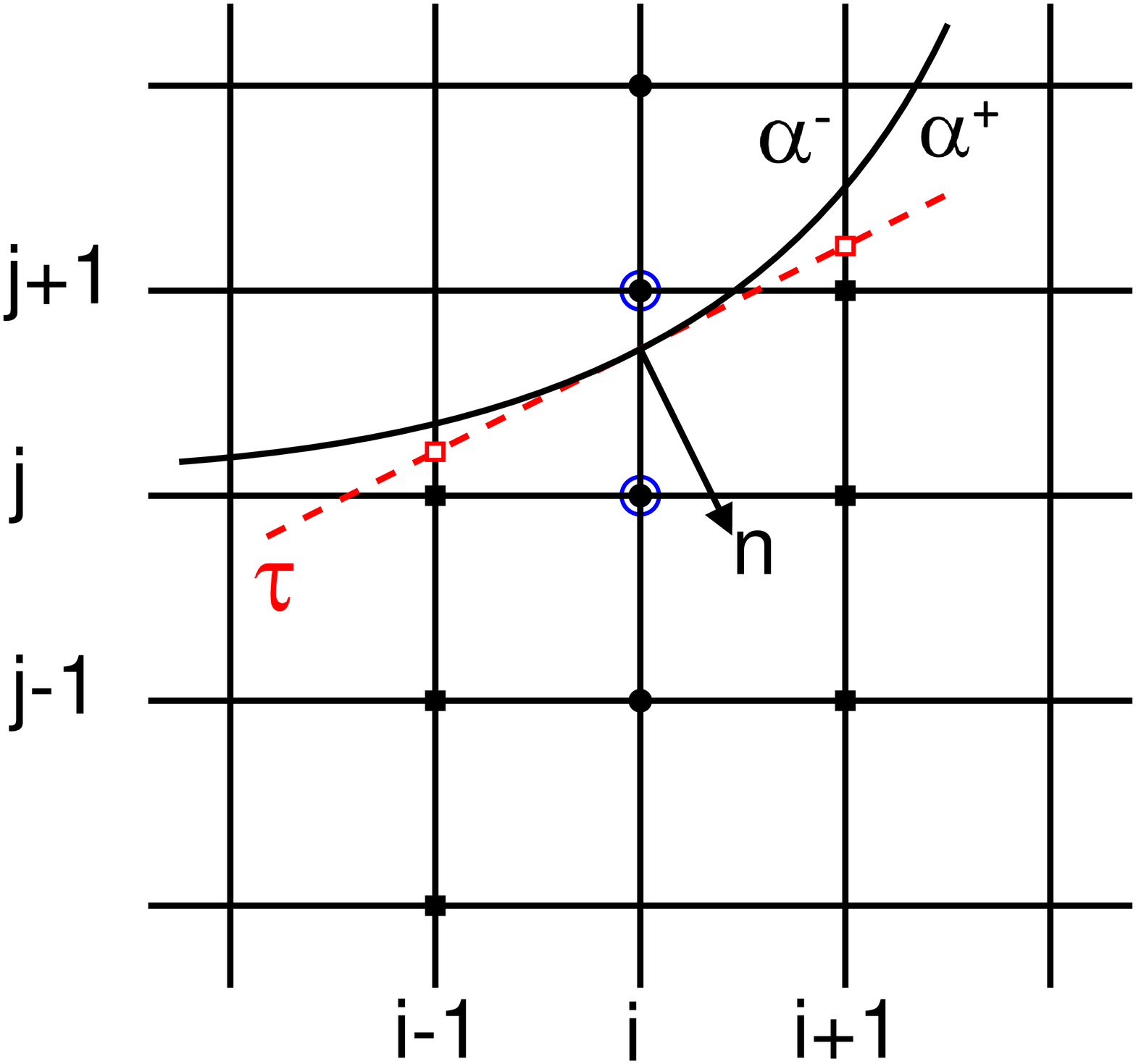,width=0.4\linewidth} & \quad \quad &
\psfig{figure=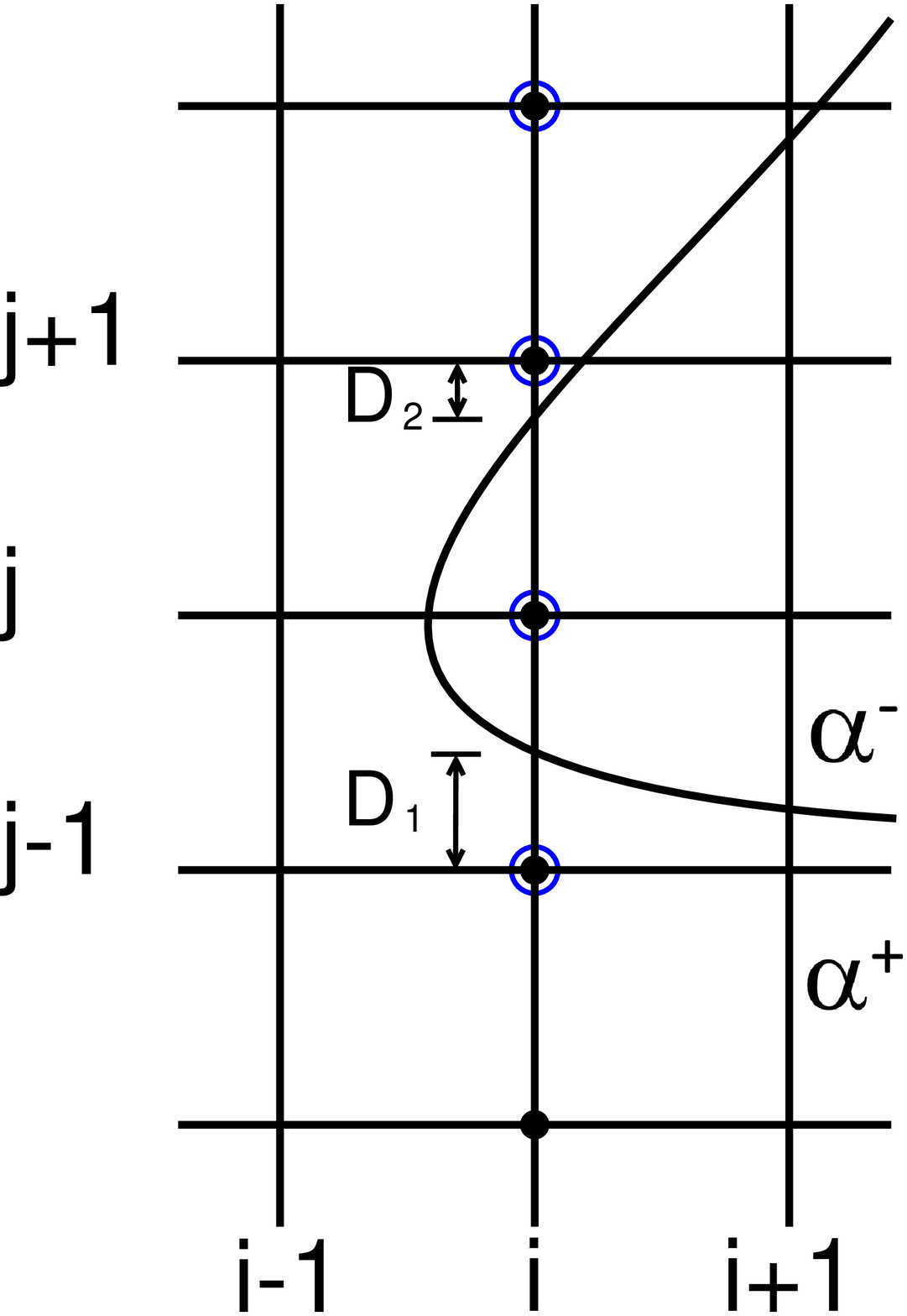,width=0.27\linewidth} \\
(a) & \quad \quad  & (b) 
\end{tabular}
\end{center}
\caption{Illustration of the MIB grid partitions.
(a). For a regular interface; (b). For a corner case. 
In both charts, the jump conditions will be discretized by using
fictitious values (open circles) and function values (filled
circles). 
In (a), the approximation of $u^+_{\tau}$ is also shown, i.e.,
it is approximated by 
two auxiliary values (open squares), then interpolated by
six function values (filled squares).
}
\label{fig.grid}
\end{figure*}

We propose a new MIB scheme to impose the decomposed 1D jump
conditions in the vicinity of the interface $\Gamma$. 
Comparing with other established interface methods, the MIB method
\cite{Zhao04,Zhou06,Zhao09} is ideally-suited to the present problem because
in the MIB, jump condition enforcement is disassociated with the derivative
discretization and can be conducted in a 1D manner along each Cartesian
direction. 
Here, we consider the MIB modification to $\delta_{yy} u^{k+1}_{i,j}$
in the ADI scheme (\ref{ADI}) as an example, 
where $(x_i,y_j)$ is an irregular node near the
interface $\Gamma$. 
A typical situation is shown in Fig. \ref{fig.grid} (a).
In the MIB scheme, to approximate function or its derivatives on one side of
interface, one never directly refers to function values from the other side. 
Instead, fictitious values from the other side of the interface will be supplied. 
Referring to Fig. \ref{fig.grid} (a), $\delta_{yy} u^{k+1}_{i,j}$ will be corrected as
\begin{equation}\label{uyy2}
\delta_{yy} u^{k+1}_{i,j}
=\frac{1}{h^2} (u^{k+1}_{i,j-1} - 2 u^{k+1}_{i,j} + \tilde{u}^{k+1}_{i,j+1}),
\end{equation}
where $\tilde{u}^{k+1}_{i,j+1}$ is a fictitious value at the node
$(x_i,y_{j+1})$. 
Similarly, $\delta_{yy} u^{k+1}_{i,j+1}$ will be modified. 
This calls for two fictitious values $\tilde{u}^{k+1}_{i,j}$ and
$\tilde{u}^{k+1}_{i,j+1}$, which will be resolved based on the jump conditions
(\ref{jump-y}).

In order to impose the jump conditions (\ref{jump-y}), we first need to
approximate $u^+_{\tau}$. 
Consider the situation shown in Fig. \ref{fig.grid} (a). 
We calculate $u^+_{\tau}$ in two steps. 
First, we calculate the intersection points between the tangential
line $\tau$ and two grid lines $x=x_{i-1}$ and $x=x_{i+1}$. 
There two auxiliary nodes are shown as open squares in Fig. \ref{fig.grid} (a).
A central difference is conducted to approximate the 
$\tau$ derivative at the interface point $(x_i,y_{\Gamma})$
using two auxiliary values of $u$.
Second, each of these two auxiliary values will be
further interpolated by using three on-grid function values, all selecting 
from $\Omega^+$.
These six nodes are shown as filled squares in Fig. \ref{fig.grid} (a). 
In this manner, $u^+_{\tau}$ is actually approximated by six grid values of $u$,
with the spatial accuracy being second order. 
In the present study, we will make use of the known values of $u$ at
the current time instant $t_k$ to estimate $u^+_{\tau}$, which avoids the 
introduction of a coupling among different $y$ grid lines at the future
time instant $t_{k+1}$. Otherwise, the 1D linear systems of the ADI algorithm
are not independent, which cannot be solved efficiently. 
We note that
the present approximation is of first order accurate in time when 
it is applied to correct $\delta_{yy} u^{k+1}_{i,j}$ and 
$\delta_{xx} u^*_{i,j}$ in (\ref{ADI}).
This is acceptable, since the temporal order of the ADI scheme (\ref{ADI})
is also one. 
In summary, the calculated $u^+_{\tau}$ will depend on six $u^{k}_{i,j}$ values 
for some nearby nodes $(x_i,y_j)$.

With the calculated $u^+_{\tau}$, 
$\hat{\psi}$ is then known at the interface point $(x_i,y_{\Gamma})$.
We will next determine $\tilde{u}^{k+1}_{i,j}$ and
$\tilde{u}^{k+1}_{i,j+1}$ by using four function values at the future time
$t_{k+1}$, i.e., 
$u^{k+1}_{i,j-1}$, $u^{k+1}_{i,j}$, $u^{k+1}_{i,j+1}$, and
$u^{k+1}_{i,j+2}$, see Fig. \ref{fig.grid} (a).
For this purpose, (\ref{jump-y}) will be discretized in the same manner of
(\ref{uyy2}), i.e., never referring to function values across the
interface $\Gamma$
\begin{align}
& w_{0,1}^+ \tilde{u}^{k+1}_{i,j} 
+ w_{0,2}^+ u^{k+1}_{i,j+1} + w_{0,3}^+ u^{k+1}_{i,j+2} \label{MIB} \\
= & w_{0,1}^- u^{k+1}_{i,j-1} + w_{0,2}^- u^{k+1}_{i,j}
+ w_{0,3}^- \tilde{u}^{k+1}_{i,j+1} + \phi,  \nonumber \\
& \alpha^+ \Big( w_{1,1}^+ \tilde{u}^{k+1}_{i,j} 
+ w_{1,2}^+ u^{k+1}_{i,j+1} + w_{1,3}^+ u^{k+1}_{i,j+2} \Big) \nonumber \\
 = & \alpha^- \Big( w_{1,1}^- u^{k+1}_{i,j-1} + w_{1,2}^- u^{k+1}_{i,j}
+ w_{1,3}^- \tilde{u}^{k+1}_{i,j+1} \Big)+ \hat{\psi}, \nonumber
\end{align}
where $w^-_{I,J}$ and $w^+_{I,J}$ for $I=0,1$ and $J=1,2,3$ are one-sided
finite difference weights, respectively, for left and right subdomains. 
Here the subscript $I$ represents interpolation ($I=0$) and the first 
derivative approximation ($I=1$), and $J$ is for grid index. 
After the discretization,  (\ref{MIB}) actually represents two algebraic equations.
By solving (\ref{MIB}), one can determine $\tilde{u}^{k+1}_{i,j}$ and
$\tilde{u}^{k+1}_{i,j+1}$ as linear combinations of 
$u^{k+1}_{i,j-1}$, $u^{k+1}_{i,j}$, $u^{k+1}_{i,j+1}$, $u^{k+1}_{i,j+2}$,
$\phi$, and $\hat{\psi}$.  
By substituting such combinations into (\ref{uyy2}) to eliminate the
fictitious value and applying the definition of $\hat{\psi}$ given in
(\ref{jump-y}), 
$\delta_{yy} u^{k+1}_{i,j}$ is now a spatially second order
accurate finite difference approximation to the double $y$ derivative,
involving four $u^{k+1}_{i,j}$ values along the $y$ direction, 
six $u^{k}_{i,j}$ values nearby, and two nonhomogeneous values $\phi$ and $\psi$.
In the present study, $\phi$ and $\psi$ will be evaluated at the time
instant $t_k$.

When one grid line intersects the interface $\Gamma$ near a rounded or 
sharp corner, the interface could be cut twice within a short distance. 
If in between these two intersection
points there is no grid node, a grid refinement is necessary. If there are at least
two grid nodes, the aforementioned matched ADI algorithm can be conducted.
However, additional corner treatments are called for the case where
only one node locates in between two intersection points, see 
Fig. \ref{fig.grid} (b). 
A MIB corner scheme is proposed to solve this problem. 
Denote two intersection points as $(x_i,y_{\Gamma 1})$ and
$(x_i,y_{\Gamma 2})$ with $y_{\Gamma 1}< y_j < y_{\Gamma 2}$. 
At  $(x_i,y_{\Gamma 1})$ and $(x_i,y_{\Gamma 2})$, 
$u^+_{\tau}$ can be calculated individually as outlined above, so that
the nonhomogeneous values $(\phi_1,\hat{\psi}_{1})$ and  
$(\phi_2,\hat{\psi}_{2})$ are known, respectively. 
The jump conditions (\ref{jump-y}) can then be imposed at 
$(x_i,y_{\Gamma 1})$ and $(x_i,y_{\Gamma 2})$ to form four algebraic equations.
Nevertheless, only three fictitious values are needed to correct 
$\delta_{yy} u$, i.e., $\tilde{u}^{k+1}_{i,j-1}$, $\tilde{u}^{k+1}_{i,j}$, and
$\tilde{u}^{k+1}_{i,j+1}$. 
This difficulty can be trivially avoided by introducing one more fictitious value. 
Referring to Fig. \ref{fig.grid} (b), we denote $D_1 = |y_{\Gamma 1} - y_{j-1}|$
and $D_2= | y_{j+1} - y_{\Gamma 2} |$. If $D_1 < D_2$, the fourth fictitious
value is chosen as $\tilde{u}^{k+1}_{i,j-2}$. Otherwise, it is selected as
$\tilde{u}^{k+1}_{i,j+2}$. With such a grid partition, each term of the 
jump conditions (\ref{jump-y}) is approximated through a third order finite
difference approximation involving four points. This ensures the overall
accuracy of the matched ADI scheme, since the solution usually undergoes
a rapid change near the corner. 
In particular, for the case shown in Fig. \ref{fig.grid} (b), the positive 
and negative terms at $(x_i,y_{\Gamma 1})$ are approximated based on 
$(u^{k+1}_{i,j-2}, u^{k+1}_{i,j-1}, \tilde{u}^{k+1}_{i,j}, u^{k+1}_{i,j+1})$ and
$(\tilde{u}^{k+1}_{i,j-1}, u^{k+1}_{i,j},\tilde{u}^{k+1}_{i,j+1},\tilde{u}^{k+1}_{i,j+2})$,
respectively, 
while those at $(x_i,y_{\Gamma 2})$ are based on 
$(u^{k+1}_{i,j-1}, \tilde{u}^{k+1}_{i,j},u^{k+1}_{i,j+1},u^{k+1}_{i,j+2})$ and
$(\tilde{u}^{k+1}_{i,j-1}, u^{k+1}_{i,j},\tilde{u}^{k+1}_{i,j+1},\tilde{u}^{k+1}_{i,j+2})$,
respectively. 
The details of the discretization are omitted here. By substituting the solved linear
combinations into (\ref{uyy2}), $\delta_{yy} u^{k+1}_{i,j}$ is approximated by
five function values 
$(u^{k+1}_{i,j-2}, u^{k+1}_{i,j-1},u^{k+1}_{i,j}, u^{k+1}_{i,j+1}, u^{k+1}_{i,j+2})$ 
and four nonhomogeneous values 
$(\phi_1,\hat{\psi}_{1}, \phi_2,\hat{\psi}_{2})$. 
Furthermore, since $u^+_{\tau}$ is evaluated at at two interface points, 
$\delta_{yy} u^{k+1}_{i,j}$ involves up to 12 nearby $u^k_{i,j}$ values,
and four jump values  $(\phi_1,\psi_{1}, \phi_2,\psi_{2})$ evaluating at time $t_k$. 

We note that the proposed MIB discretization for the regular interface case and
the corner case need to be conducted only once
at the beginning of the simulation, because the geometric domain, grid, and
finite difference approximations of the jump conditions are all time invariant.
In fact, the entries of the discrete operators $\delta_{xx} u$ and
$\delta_{yy} u$ are all time independent and can be pre-determined. 
At each time step, one just needs to update nonhomogeneous values
$\phi$, $\bar{\psi}$ and $\hat{\psi}$ for time dependent jump conditions. 

\subsection{Fast algebraic solution}
\begin{figure*}[!tb]
\begin{center}
\begin{tabular}{ccc}
\psfig{figure=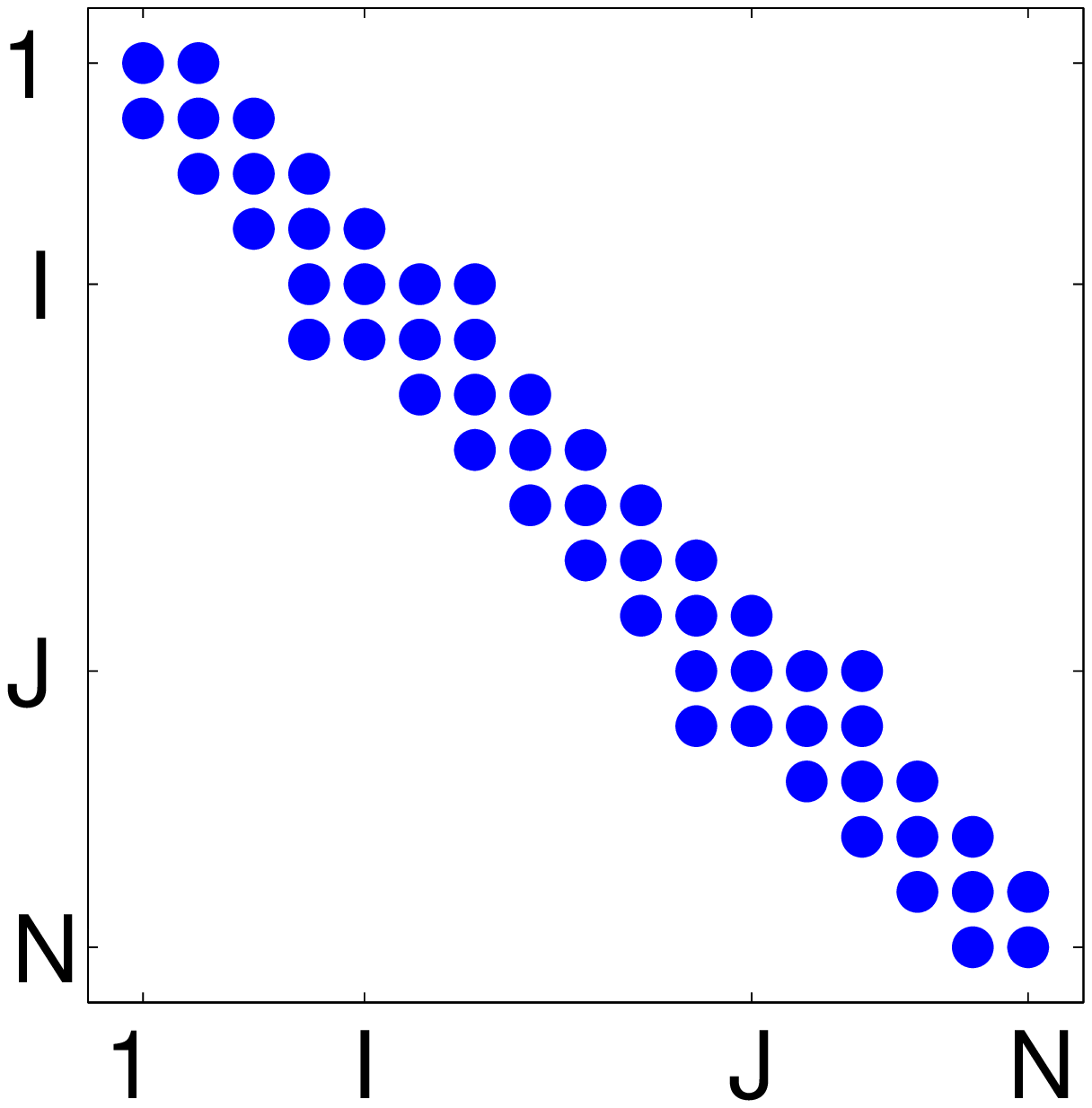,width=0.3\linewidth} & \quad \quad &
\psfig{figure=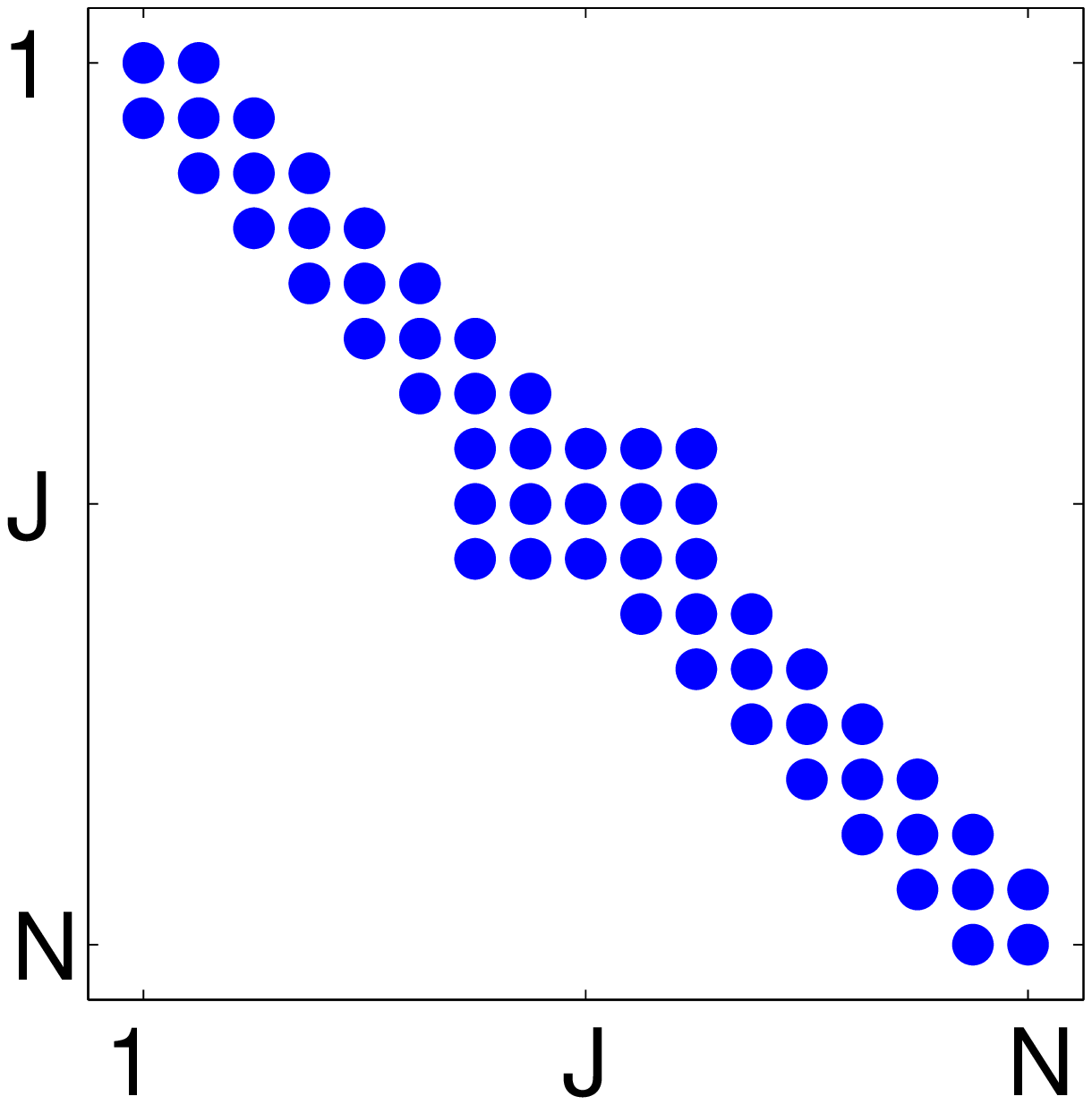,width=0.3\linewidth} \\
(a) & \quad \quad  & (b) 
\end{tabular}
\end{center}
\caption{Matrix structures of the matched ADI algorithm.
(a). For a regular interface; (b). For a corner case. 
}
\label{fig.matrix}
\end{figure*}

In the proposed matched ADI algorithm, the 1D linear systems underlying 
(\ref{ADI}) are actually independent from each other, 
so that we can solve them separately. 
Without the loss of generality, we denote the 1D linear system to be solved 
in one ADI step as
\begin{equation}\label{Axb}
{\bf A} {\bf x} = {\bf b},
\end{equation}
where ${\bf x}$ represents unknown $u^{k+1}$ or $u^*$ values on one 
$x$ or $y$ grid line. The matrix ${\bf A}$ is of dimension $N$ by $N$, where 
$N$ could be $N=N_x$ or $N=N_y$.
The vector ${\bf b}$ contains all the right hand side terms. 
Nevertheless, ${\bf A}$ will be  non-tridiagonal after the MIB treatment
so that new algebraic solvers have to be developed to maintain the 
overall computational efficiency.
In particular, for a regular interface case, 
one grid line cuts the interface $\Gamma$ twice at a well separated distance. 
There are totally four irregular nodes sandwiched the interface.
For each of them, the actual band width becomes four, 
while two nonhomogeneous values $\phi$ and $\psi$ shall be added into ${\bf b}$. 
A typical matrix structure of ${\bf A}$ is shown in Fig. \ref{fig.matrix} (a).
For the corner interface case, a similar analysis indicates that
the band width of three consecutive irregular nodes
are changed from three to five, see  Fig. \ref{fig.matrix} (b).

Since the change of the band-structure is not too dramatic, 
the linear system (\ref{Axb}) could be
solved by the Woodbury formula \cite{NR}. 
We consider the regular interface case as an example.  
Denote the indices of four irregular nodes to be $I$, $I+1$, $J$, and $J+1$. 
See Fig. \ref{fig.matrix} (a).
The extra coefficients of ${\bf A}$ can be accounted for by defining two 
$N \times 2$ matrices ${\bf P}$ and ${\bf Q}$ with four nonzero elements each:
\begin{align*}
& {\bf P}_{I,1}=1, \quad {\bf P}_{I+1,1}=1,  \quad {\bf P}_{J,2}=1, 
\quad {\bf P}_{J+1,2}=1, \\
& {\bf Q}_{I-1,1}={\bf A}_{I+1,I-1},  \quad {\bf Q}_{I+2,1}={\bf A}_{I,I+2}, \quad
{\bf Q}_{J-1,2}={\bf A}_{J+1,J-1},  \quad {\bf Q}_{J+2,2}={\bf A}_{J,J+2}.
\end{align*}
We then have ${\bf A}={\bf T} + {\bf P}{\bf Q}^T$, where ${\bf T}$
is a tridiagonal matrix. Thus, we have analytically
\begin{equation}
{\bf A}^{-1}=({\bf T} + {\bf P}{\bf Q}^T)^{-1}
={\bf T}^{-1} - [ {\bf T}^{-1} {\bf P} 
({\bf 1} +  {\bf Q}^T {\bf T}^{-1} {\bf P})^{-1} {\bf Q}^T {\bf T}^{-1}],
\end{equation}
by the Woodbury formula \cite{NR}.
In other words, the inversion of ${\bf A}$ can be carried out through applying the
Thomas algorithm three times by solving some auxiliary systems about ${\bf T}$.

However, the Woodbury formula is not used in our computation, 
because a more efficient algebraic procedure is available. 
For example, for the regular interface case, 
four elementary row operations are simply conducted. 
In particular, by denoting the $I^{\rm th}$ row of ${\bf A}$ as ${\bf R}_I$,
the elementary row operation 
${\bf R}_{I+1} - {\bf A}_{I+1,I-1}/{\bf A}_{I,I-1} {\bf R}_{I}$
will zero the entry at the position $(I+1,I-1)$. 
Similar row operations are conducted to vanish other three positions, so that
a tridiagonal system is formed. 
In a similar manner, the matrix ${\bf A}$ of the corner interface can be
treated by six row operations to reduce to a  tridiagonal one. 
After the Gauss elimination, 
the Thomas algorithm \cite{ADIbook}
is applied only once to solve the adjusted system. 

The proposed matched ADI algorithm is very efficient. 
To solve the 1D linear system within each inner ADI step, 
the flop counts are essentially due to the Thomas algorithm, 
i.e., on the order of $O(N)$,
because the overhead for the Gauss elimination is very small and does
not grow with $N$. 
Thus, like the standard 2D ADI algorithm, the complexity 
of advancing one time step in the proposed matched ADI schemes is about
$O(N^2)$.
Moreover, due to the excellent stability of the matched ADI scheme,
one can simply fix $\Delta t$ to be  on the order of $h$. 
Consequently, the complexity of
entire time integration will be on the order of $O(N^3)$ for solving
 parabolic interface problems with $N^2$ unknowns.

\subsection{Stability analysis}\label{sec.stab}
As discussed above that the Douglas ADI scheme (\ref{ADI}) is a higher
order perturbation of the implicit Euler scheme (\ref{imEuler}). 
In the present subsection, we will analyze the stability of the 
implicit Euler scheme with the proposed MIB spatial discretization, 
because in this case, it is relatively easier 
to construct 2D matrices for the spectrum analysis. 
The stability of the matched ADI scheme is essentially determined by 
that of the matched Euler scheme.

We first establish a vector notation
for the proposed MIB spatial discretization.
Denote ${\bf U}^{k} = [ u^k_{1,1}, u^k_{2,1}, \ldots, u^k_{N_x,1}, 
u^k_{1,2}, u^k_{2,2}, \ldots, u^k_{N_x,2}, \ldots]^T$, 
which is a vector of the length $N_x \times N_y$, 
containing all $u$ values at the time $t_k$.  
The second order $x$ derivative of ${\bf U}^{k+1}$ can be expressed as
\begin{equation}\label{Mat_x}
\frac{\partial ^2}{\partial x^2} {\bf U}^{k+1} \approx
{\bf D}_{xx} {\bf U}^{k+1} + {\bf \bar{B}} {\bf U}^{k} + {\bf \bar{\Phi}}^k, 
\end{equation}
where ${\bf D}_{xx} $ and ${\bf \bar{B}}$ are matrices of the dimension 
$N_x \times N_y$ by $N_x \times N_y$, while ${\bf \bar{\Phi}}^k$ is a vector
of the length $N_x \times N_y$.  
The matrix ${\bf D}_{xx} $ is a perturbation of the standard matrix for the central
difference approximation, while ${\bf \bar{B}}$  is due to the approximation of 
$u^+_{\tau}$ by some $u^k_{i,j}$ values. 
Here ${\bf \bar{\Phi}}^k$  is a correction term, 
based on linear combinations of the nonhomogeneous values $\phi$ and $\psi$. 
For a regular node $(x_i,y_j)$, the corresponding rows of ${\bf \bar{B}}$ 
and ${\bf \bar{\Phi}}^k$  have only zero entries, while that of ${\bf D}_{xx} $ has three non-zero
entries centered at the diagonal, i.e., $1/h^2$, $-2/h^2$ and $1/h^2$. 
For an irregular node in the non-corner case, the corresponding row of
${\bf D}_{xx} $ and ${\bf \bar{B}}$ has four and six non-zero elements, respectively, 
whereas in a corner case, there are five and twelve non-zero elements, respectively,
for ${\bf D}_{xx} $ and ${\bf \bar{B}}$ . 
We note that sparse elements of ${\bf \bar{B}}$ are distributed in a rather random
fashion -- their locations depend on the interface geometry and grid size. Our MIB
code will automatically calculate them. 
However, the sparse structure of ${\bf D}_{xx} $  can be well predicted. Essentially,
${\bf D}_{xx} $  has $N_y$ nonzero blocks along the diagonal. 
When the grid line $x=x_i$ does not cut the interface $\Gamma$, 
the corresponding diagonal block of ${\bf D}_{xx} $ is simply a tridiagonal sub-matrix.
Otherwise,
the corresponding diagonal block  will take the form showing in either
Fig. \ref{fig.matrix} (a) or (b), depending on whether this is a corner case or not.
Similarly, the second order $y$ derivative of ${\bf U}^{k+1}$ is approximated as
\begin{equation}\label{Mat_y}
\frac{\partial ^2}{\partial y^2} {\bf U}^{k+1} \approx
{\bf D}_{yy} {\bf U}^{k+1} + {\bf \hat{B}} {\bf U}^{k} + {\bf \hat{\Phi}}^k, 
\end{equation}
with ${\bf D}_{yy}$, ${\bf \hat{B}}$, and $ {\bf \hat{\Phi}}^k$ being 
appropriately defined.

The vector form of the implicit Euler scheme (\ref{imEuler}) can then be given as
\begin{equation}\label{stab1}
\left( \frac{1}{\alpha} {\bf I} - \Delta t {\bf D}_{xx}- \Delta t {\bf D}_{yy} \right) 
{\bf U}^{k+1}  = \left( \frac{1}{\alpha} {\bf I}  + \Delta t {\bf \bar{B}} 
+ \Delta t {\bf \hat{B}} \right) {\bf U}^{k} + \Delta t {\bf \bar{\Phi}}^k 
+ \Delta t {\bf \hat{\Phi}}^k + \frac{\Delta t }{\alpha}  {\bf F}^{k+1},
\end{equation}
where ${\bf I}$ is the identity matrix and ${\bf F}^{k+1}$ represents the 
source term. We can rewrite (\ref{stab1}) into a more compact form
\begin{equation}\label{stab2}
 {\bf D} {\bf U}^{k+1} = {\bf B} {\bf U}^{k} + {\bf C},
\end{equation}
where ${\bf D} =\frac{1}{\alpha} {\bf I} - \Delta t {\bf D}_{xx}- \Delta t {\bf D}_{yy}$,
${\bf B}=\frac{1}{\alpha} {\bf I}  + \Delta t {\bf \bar{B}} + \Delta t {\bf \hat{B}}$, and
${\bf C}=\Delta t {\bf \bar{\Phi}}^k + 
\Delta t {\bf \hat{\Phi}}^k + \frac{\Delta t }{\alpha}  {\bf F}^{k+1}$. 
Recall again that ${\bf D}$ is slightly modified from the standard matrix for the central
finite difference, while ${\bf B}$ is resulting from the approximation of $u^+_{\tau}$
by $u^k_{i,j}$ values at various interface points. 
By taking an inverse, (\ref{stab2}) becomes
\begin{equation}\label{stab3}
{\bf U}^{k+1} = {\bf D}^{-1}  {\bf B} {\bf U}^{k} + {\bf D}^{-1} {\bf C}
={\bf M} {\bf U}^{k} + {\bf D}^{-1} {\bf C}.
\end{equation}
Thus, the stability of the entire spatial-temporal discretization depends on the 
magnifying matrix ${\bf M} = {\bf D}^{-1}  {\bf B}$. 
In particular, it depends on the spectral radius $\rho$ of ${\bf M}$, which
is defined to be
\begin{equation}
\rho ({\bf M}) = \max_j | \lambda_j |,
\end{equation}
where $ \lambda_j$ are eigenvalues of ${\bf M}$.

In general, an analytical spectrum analysis of ${\bf M}$ is extremely difficult,
because the finite difference formulas underlying ${\bf D}$ and ${\bf B}$
depend on the positions of grid nodes and interface intersection points. Thus,
in the present study, the leading eigenvalues of ${\bf M}$ will be calculated
numerically. This enables us to directly examine the impact of
parameters, such as  $h$, $\Delta t$ and $\alpha$, on the stability. 
Moreover, we will also investigate how a complex geometry could affect the
stability.

The stability of the MIB spatial discretization combined with explicit time
integrations has been analyzed in \cite{Zhao04,Zhao07,Zhao09}
for solving both interface and boundary closure problems. 
It is known that an instability could occur, if a too asymmetric MIB 
finite difference approximation is involved \cite{Zhao04,Zhao09},  
while a severely asymmetric approximation will produce
spurious modes in the spectrum analysis \cite{Zhao07}. 
When a complicated interface $\Gamma$ is studied, 
a severely asymmetric approximation could be encountered, if 
$\Gamma$ and/or its tangential lines cut a grid line at
a point that is very close to a node. 
However, one cannot predict when this will happen,
because $\Gamma$ will intersect the grid in a random manner.
To evaluate the impact of such a geometrical effect on the stability, 
we will consider a fixed interface and test various 
different mesh sizes in the next section.
This allows us to see if the complex geometry will compromise 
the stability of the present implicit time stepping method.

\section{Numerical experiments}\label{Sec:Numerical}
In this section, we investigate the stability and accuracy of
the proposed matched ADI algorithm for
solving 2D parabolic interface problems with different jump conditions
and interface geometries. 
Piecewisely defined analytical solutions will be constructed  in each example. 
The initial solution is chosen according the analytical solution at $t=0$.
The ADI time stepping will be carried out until a stopping time $t=T$. 
Without the loss of generality, a square domain $[-D,D]\times[-D,D]$ 
with the Dirichlet boundary condition is considered in all examples. 
Here, the boundary data is simply given by the analytical solutions. 
Similarly, the jump conditions defining function and flux jumps
across the interfaces are also calculated according to the given
analytical solutions. 
For simplicity, the mesh sizes in both $x$ and $y$ directions are chosen 
to be the same, i.e., $N=N_x=N_y$ with $h=\Delta x =\Delta y$. 
The domain size $D$ is usually chosen as a non-integer so that
the corner interface case could be encountered in a coarse grid. 
This enables us to fully validate the proposed matched ADI algorithm. 
Numerical errors in $L_\infty$ and $L_2$ norms are
reported in all examples. 

In all examples, the proposed matched ADI method is found to be 
unconditionally stable, through both direct numerical verifications and 
eigenvalue stability analysis. To save the space, the detailed stability
analysis will be presented only for the last numerical example,
in which the most complicated interface geometry and the most complicated
jump conditions will be studied. 

{\bf Example 1.}  We first study a circular interface problem with a
continuous solution. 
Consider a square domain $[-D,D]\times[-D,D]$ 
with a circular interface $r^2=x^2+y^2=1$. 
The piecewise coefficient is defined to be $\alpha^-=1$ and $\alpha^+=10$,
respectively, for $r<1$ and $r \ge 1$. 
The analytical solution to the heat equation is designed to be
\begin{equation*}
u(x,y,t)= \begin{cases}
\Big( \frac{r^6 -1}{\alpha^-} - \frac{3}{\alpha^+} \Big) \cos(t),
\quad      \mbox{ if }   r <1 \\
-\frac{3}{\alpha^+ r^2} \cos(t),
\quad \quad \quad \quad   \mbox{ if } r \ge 1,
\end{cases}
\end{equation*}
so the the jump conditions are simply $[u]=0$ and $[\alpha u_n]=0$. 
The source term is then given as
\begin{equation*}
f(x,y,t)= \begin{cases}
-\Big( \frac{r^6 -1}{\alpha^-} - \frac{3}{\alpha^+} \Big) \sin(t) 
-36 r^4 \cos(t),
\quad      \mbox{ if }   r <1 \\
\frac{3}{\alpha^+ r^2} \sin(t) +12 r^{-4} \cos(t),
\quad \quad \quad \quad   \quad \mbox{ if } r \ge 1.
\end{cases}
\end{equation*}
In this example, the domain size is set to be $D=1.99$ and
the stopping time is fixed as $T=2$.

\begin{figure*}[!tb]
\begin{center}
\begin{tabular}{cc}
\psfig{figure=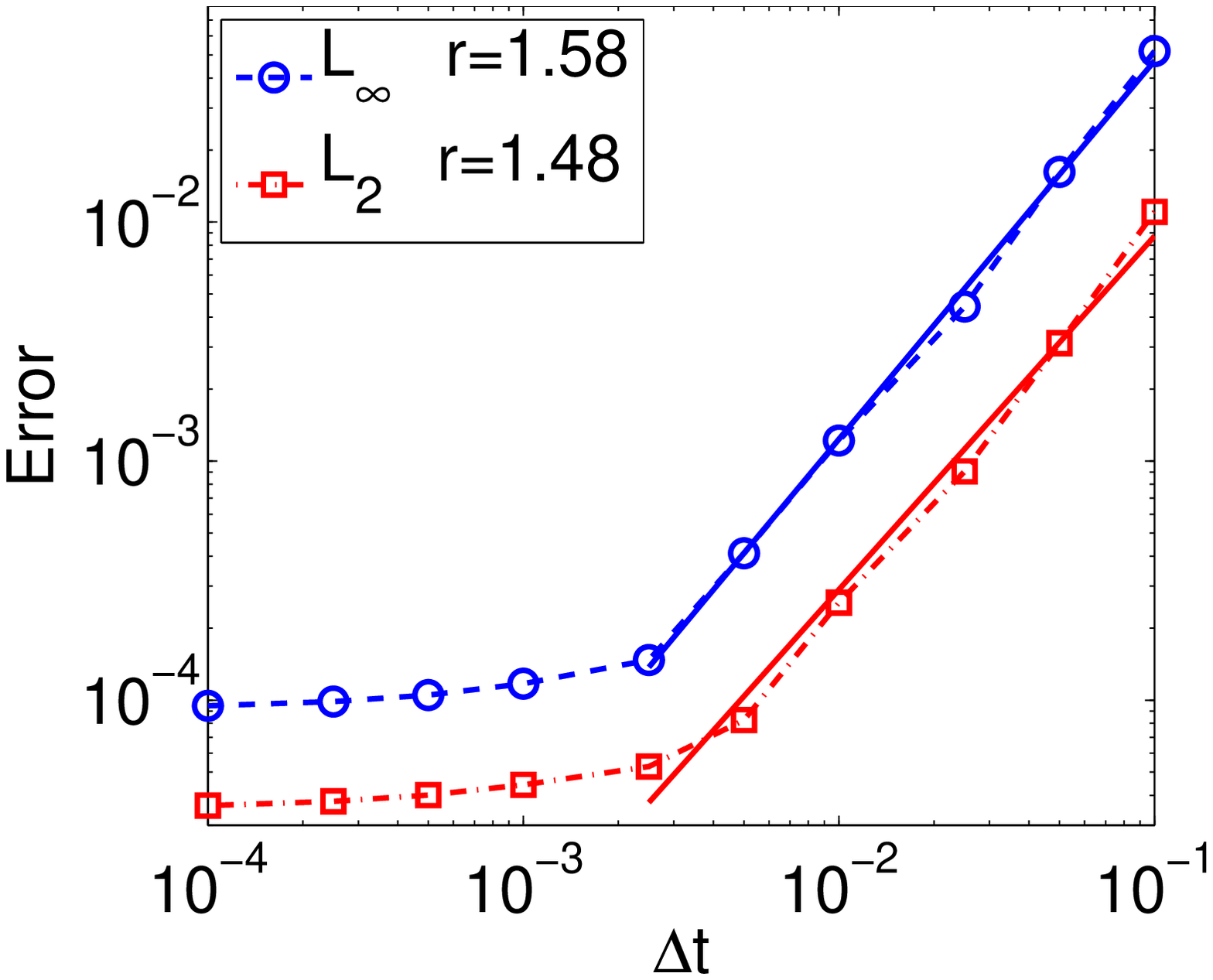,width=0.45\linewidth} & 
\psfig{figure=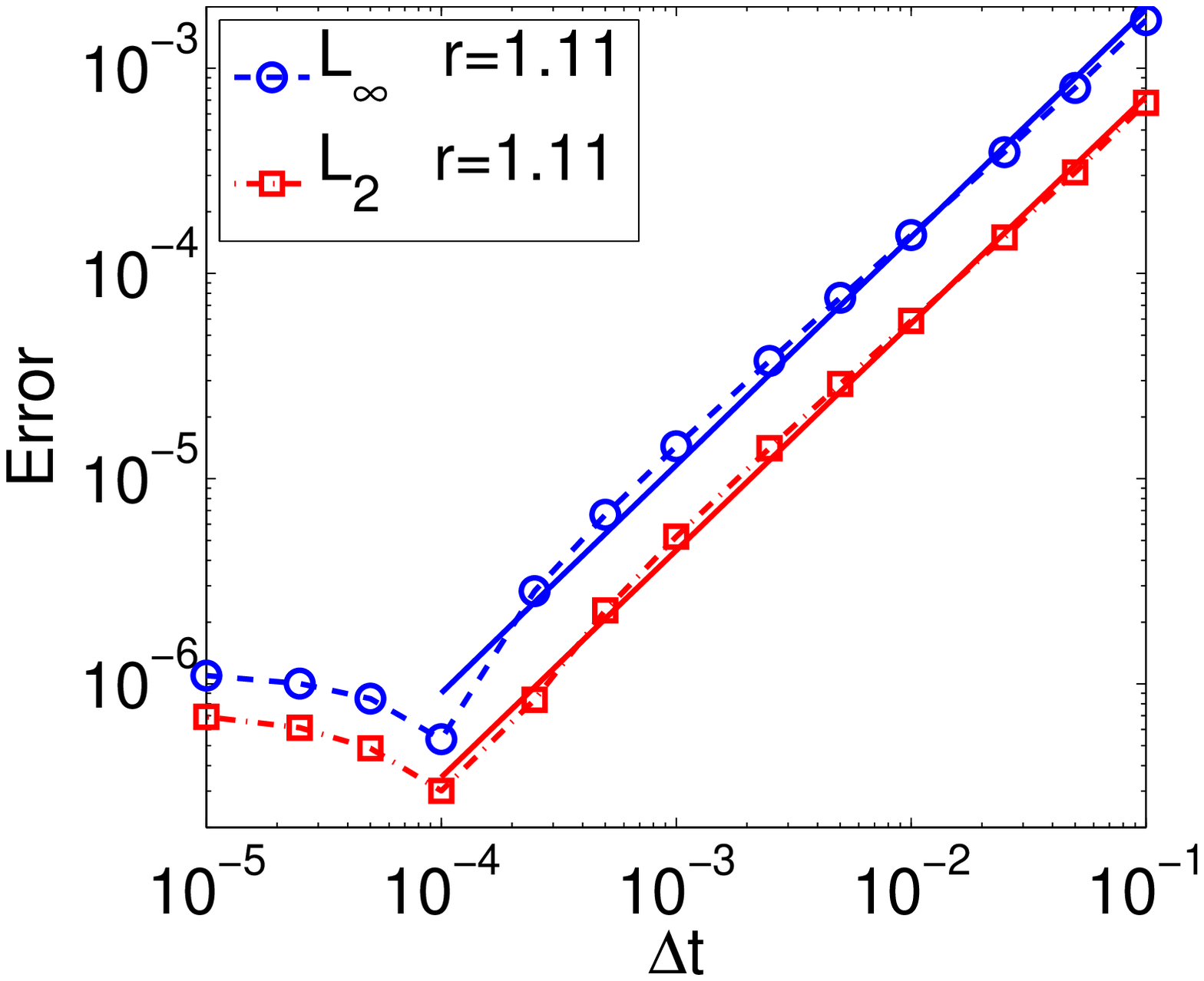,width=0.45\linewidth} \\
(a) &  (b) 
\end{tabular}
\end{center}
\caption{Temporal convergence tests. (a).  Example 1;
(b). Example 2. 
}
\label{fig.ex12}
\end{figure*}

We first examine the temporal convergence. 
The proposed matched ADI algorithm is found to be unconditionally stable for
all tested $h$ and $\Delta t$.  
By taking $N=321$, 
the numerical errors generated by using different $\Delta t$ are shown in 
Fig. \ref{fig.ex12} (a). 
A similar pattern can be observed for both $L_\infty$ and $L_2$ errors.
i.e., the temporal errors become smaller and smaller
until they are limited by the accuracies of the spatial discretization. 
The temporal convergence order of the matched ADI algorithm
can be analyzed via considering 
errors before reaching the limiting precision. 
For these errors, a linear least-squares fitting \cite{Zhao09} is conducted
in the log-log scale. The fitted convergence lines are shown as solid lines in 
Fig. \ref{fig.ex12} (a). Moreover, the fitted slope essentially represents
the numerical convergence rate $r$ of the scheme. 
The temporal order in $L_\infty$ and $L_2$ norms is found to be
$r=1.58$ and $r=1.48$, respectively. 
In other word, the numerically detected temporal oder of the matched ADI algorithm
is about half order higher than its theoretical design. 
This is perhaps because the present solution is continuous.
The superconvergence  of the Douglas ADI scheme for smooth solutions has also
been observed in other literature studies \cite{Geng13,Zhao14}. 

\begin{table}[!t]
\caption{Spatial convergence tests for first two examples. }
\label{table.ex12}
\begin{center}
\begin{tabular}{|l|l|l|l|l|l|l|l|l|}
\hline
& \multicolumn{4}{c|}{Example 1} & \multicolumn{4}{c|}{Example 2} \\
\cline{2-5} \cline{6-9}
& \multicolumn{2}{c|}{$L_\infty$} & \multicolumn{2}{c|}{$L^2$} 
& \multicolumn{2}{c|}{$L_\infty$} & \multicolumn{2}{c|}{$L^2$} \\
\cline{2-3} \cline{4-5} \cline{6-7} \cline{8-9}
$N$ & error & order & error & order & error & order & error & order  \\
\hline
21   &  1.92E-2 &         &  6.75E-3 &         & 1.78E-4 &         & 9.15E-5 & \\
41   &  4.49E-3 & 2.09 &  1.51E-3 & 2.16 & 4.38E-5 & 2.02 & 2.76E-5 & 1.73 \\
81   &  1.34E-3 & 1.74 &  4.97E-4 & 1.61 & 1.77E-5 & 1.31 & 1.13E-5 & 1.28 \\
161 &  3.78E-4 & 1.83 &  1.46E-4 & 1.77 & 3.07E-6 & 2.53 & 2.08E-6 & 2.45 \\
321 &  9.47E-5 & 2.00 &  3.62E-5 & 2.01 & 5.38E-7 & 2.51 & 3.01E-7 & 2.79\\
\hline
\end{tabular}
\end{center}
\end{table}

We next quantitatively examine the spatial accuracy. 
A small enough $\Delta t= 10^{-4}$ is employed so that the temporal error
can be neglected in the present study.
The numerical errors of the matched ADI method
for different mesh size $N$ are listed in Table \ref{table.ex12}. 
Based on successive mesh refinements,
the numerically calculated convergence rates are also reported for both
error measurements. 
It can be seen that the matched ADI algorithm achieves the second order of 
accuracy in both $L_\infty$ and $L_2$ norms 
for the present parabolic interface problem.

{\bf Example 2.} 
We next consider a circular interface problem with constant jump values. 
The interface $\Gamma$ is defined as $r=0.5$ and the domain size is
set as $D=0.99$. 
The diffusion coefficient is chosen as $\alpha^- =2$ and $\alpha^+=10$,
respectively, for $r < 0.5$ and $r \ge 0.5$. 
The source term is defined to be 
\begin{equation*}
f(x,y,t)= \begin{cases}
-\sin(t) - 4 \alpha^-,  \quad  \quad \quad    \mbox{ if }   r <0.5 \\
-\sin(t) - 8 r^2 -4, 
\quad \quad   \mbox{ if } r \ge 0.5.
\end{cases}
\end{equation*}
The analytical solution can then be given as 
\begin{equation*}
u(x,y,t)= \begin{cases}
 \cos(t) + r^2 -1,  \quad  \quad  \quad  \quad  \quad
 \quad  \quad    \quad  \quad   \quad     \mbox{ if }   r <0.5 \\
\cos(t) + \frac{1}{4}(1- \frac{9}{8 \alpha^+}) + 
\frac{1}{\alpha^+} (\frac{r^4}{2} + r^2),
\quad   \mbox{ if } r \ge 0.5.
\end{cases}
\end{equation*}
It can be verified that the jump values are constants along the interface $\Gamma$
and are time invariant. 
In particular, we have $[u]=1$ and $[\alpha u_n]=-0.75$. 

For all tested $h$ and $\Delta t$, the matched ADI method is again found to
be unconditionally stable for this example. 
By choosing $N=321$ and $T=1$, 
the temporal accuracies are analyzed in 
Fig. \ref{fig.ex12} (b). It can be seen that both $L_\infty$ and $L_2$ errors
decrease uniformly until the limiting precisions of the spatial discretization are
reached. After that, by using a smaller $\Delta t$, the error become slightly
larger. The least-squares error analysis is also conducted. With a rate of $1.11$
for both error norms, the matched ADI method clearly attains the first order of
accuracy in time for this example with a discontinuous solution. 
By using a sufficiently small $\Delta t=10^{-4}$, the spatial accuracies of the
matched ADI method are investigated in Table \ref{table.ex12}. 
The spatial convergence is not very uniform for this example. 
However, the overall order of the matched ADI method is still around two
in both error measurements. 

{\bf Example 3.} 
To further explore the potential of the proposed ADI method, we consider
a circular interface problem with general jump values. 
The interface is also defined as $r=0.5$ with domain size $D=0.99$. 
By taking $\alpha^- = 1$ and $\alpha^+=10$, the analytical solution is chosen as
\begin{equation*}
u(x,y,t)= \begin{cases}
 \cos(t) + \exp(x^2+y^2),  \quad  \quad   \mbox{ if }   r <0.5 \\
\cos(t) + \sin(kx) \cos(ky),
\quad   \mbox{ if } r \ge 0.5,
\end{cases}
\end{equation*}
where the wavenumber is chosen as $k=2$. 
The source term is given as
\begin{equation*}
f(x,y,t)= \begin{cases}
-\sin(t) -  4 \alpha^- \exp(x^2+y^2)(x^2+y^2+1),  \quad    \mbox{ if }   r <0.5 \\
-\sin(t) + 2 \alpha^+ k^2 \sin(kx)\cos(ky), 
\quad \quad  \quad \quad  \quad \mbox{ if } r \ge 0.5.
\end{cases}
\end{equation*}
The jump conditions at a interface point 
$(x,y)=(\frac{1}{2} \cos \theta, \frac{1}{2} \sin \theta)$ 
can be derived from the analytical solution
\begin{align*}
[u] =& \sin (\frac{k}{2} \cos \theta) \cos (\frac{k}{2} \sin \theta) - \exp(\frac{1}{4}), \\
[\alpha u_n] =&
 \alpha^+ k \cos \theta \cos (\frac{k}{2} \cos \theta) \cos (\frac{k}{2} \sin \theta) \\
- & \alpha^+ k \sin \theta \sin (\frac{k}{2} \cos \theta) \sin (\frac{k}{2} \sin \theta)
-\alpha^- \exp(\frac{1}{4}), \\
[u_{\tau}] =&  
- k \sin \theta \cos (\frac{k}{2} \cos \theta) \cos (\frac{k}{2} \sin \theta)  
- k \cos \theta \sin (\frac{k}{2} \cos \theta) \sin (\frac{k}{2} \sin \theta) ,
\end{align*}
where the third jump condition $[u_{\tau}]=\phi_{\tau}$ is derived from the function 
jump $[u]=\phi$. 
We note that the present jump conditions are quite general in the sense that 
all jump values are functions of space, even though they are time independent.

\begin{figure*}[!tb]
\begin{center}
\begin{tabular}{cc}
\psfig{figure=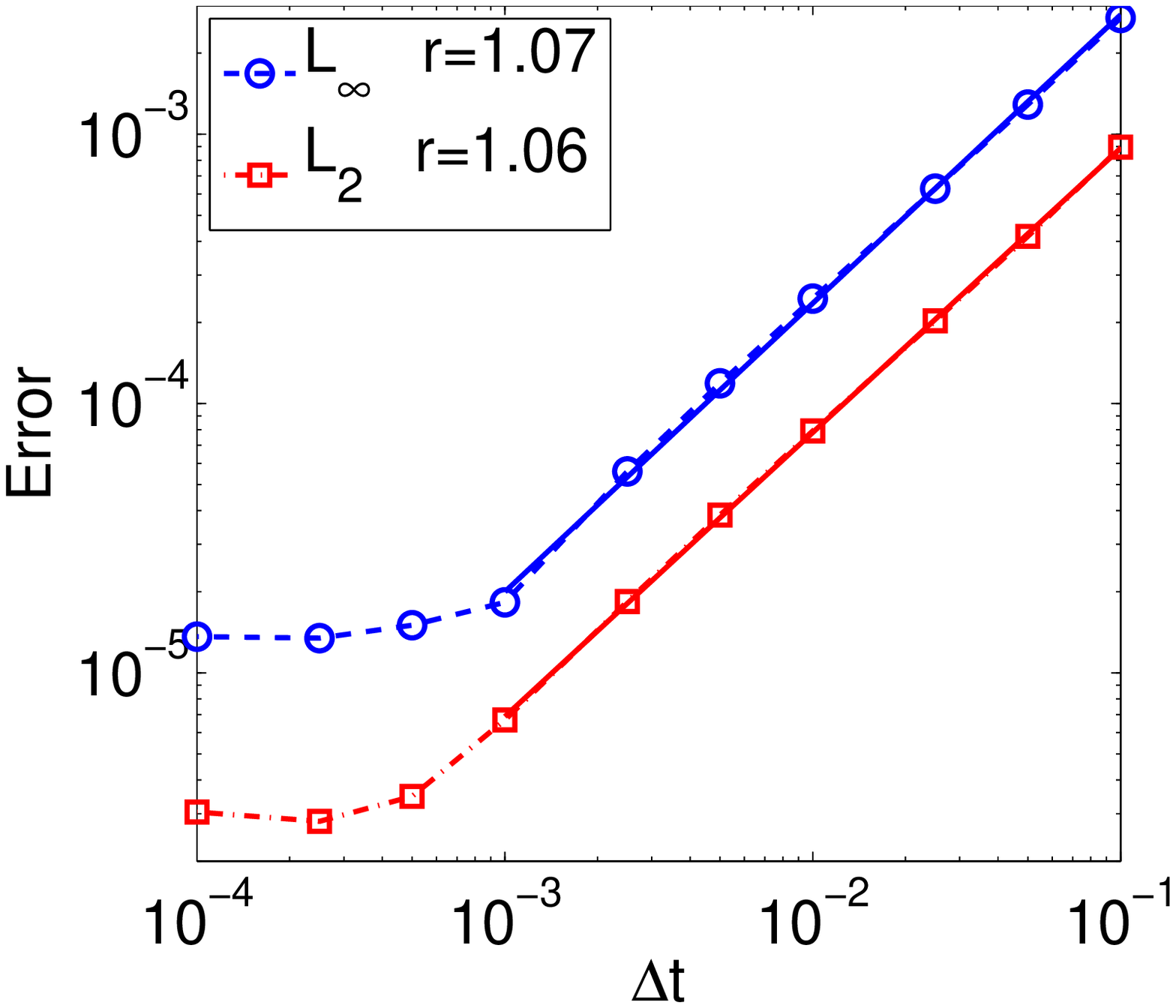,width=0.45\linewidth} & 
\psfig{figure=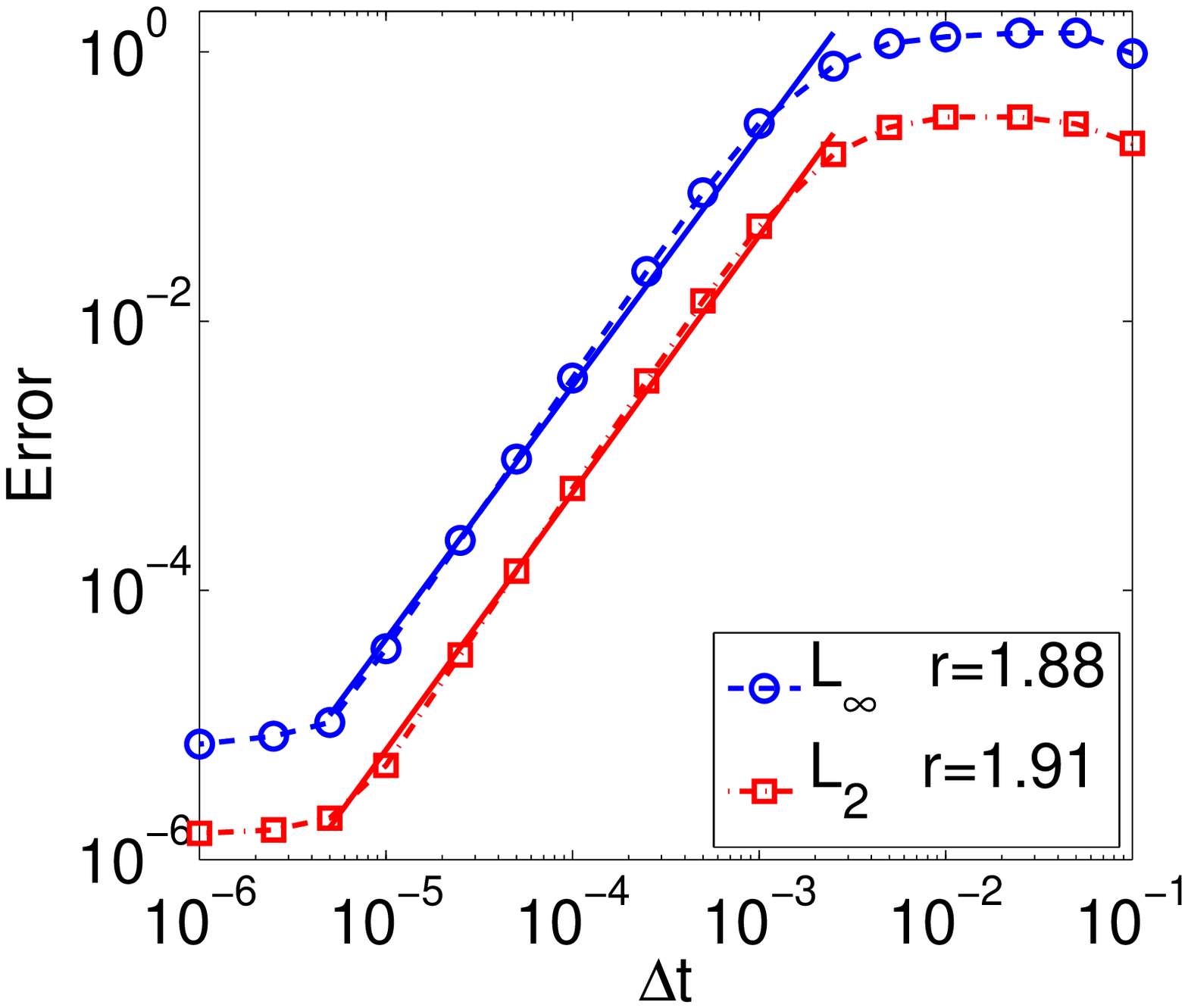,width=0.45\linewidth} \\
(a) &  (b) 
\end{tabular}
\end{center}
\caption{Temporal convergence tests. (a).  Example 3;
(b). Example 4. 
}
\label{fig.ex34}
\end{figure*}

\begin{figure*}[!tb]
\begin{center}
\begin{tabular}{cc}
\psfig{figure=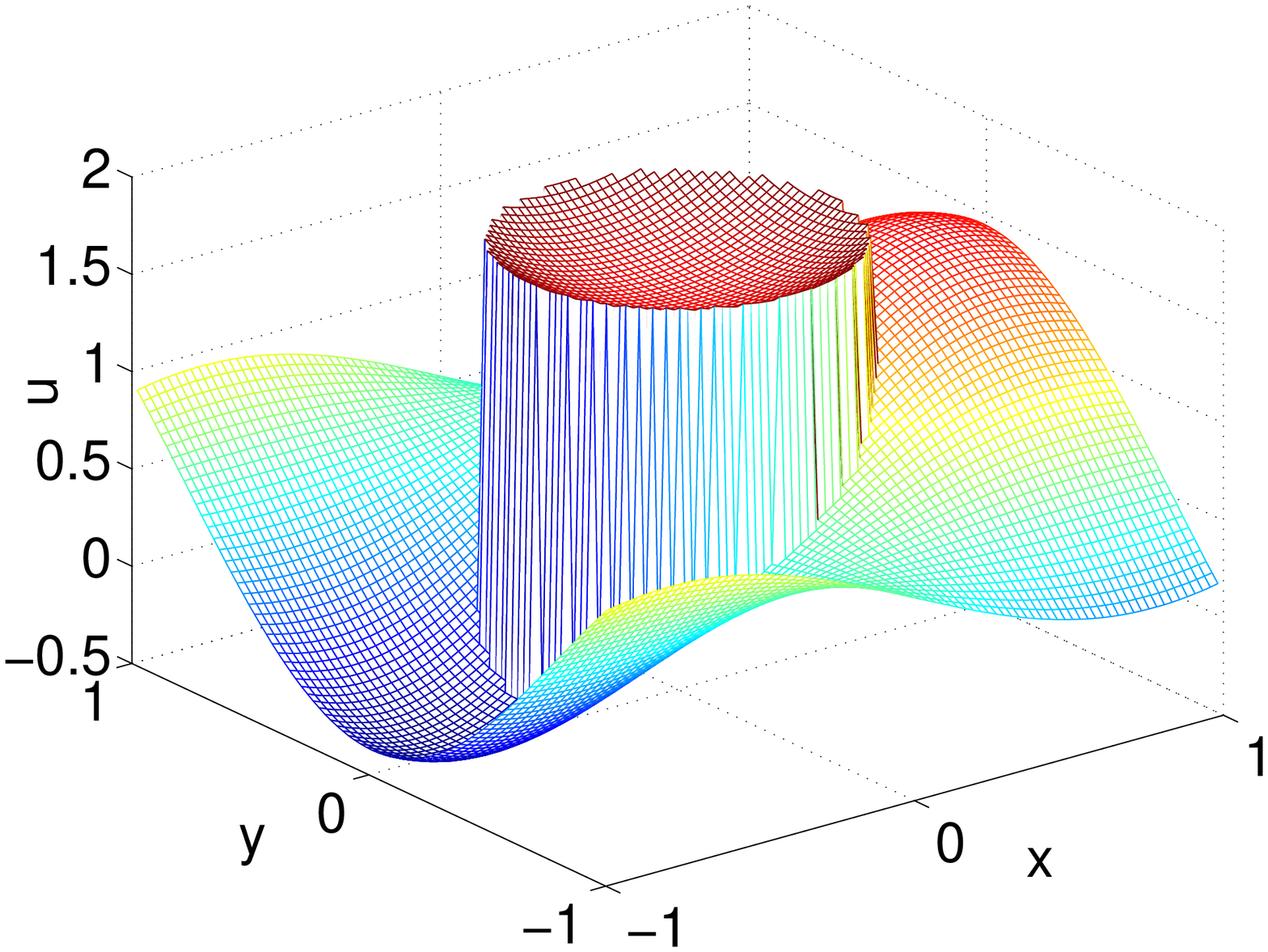,width=0.45\linewidth} & 
\psfig{figure=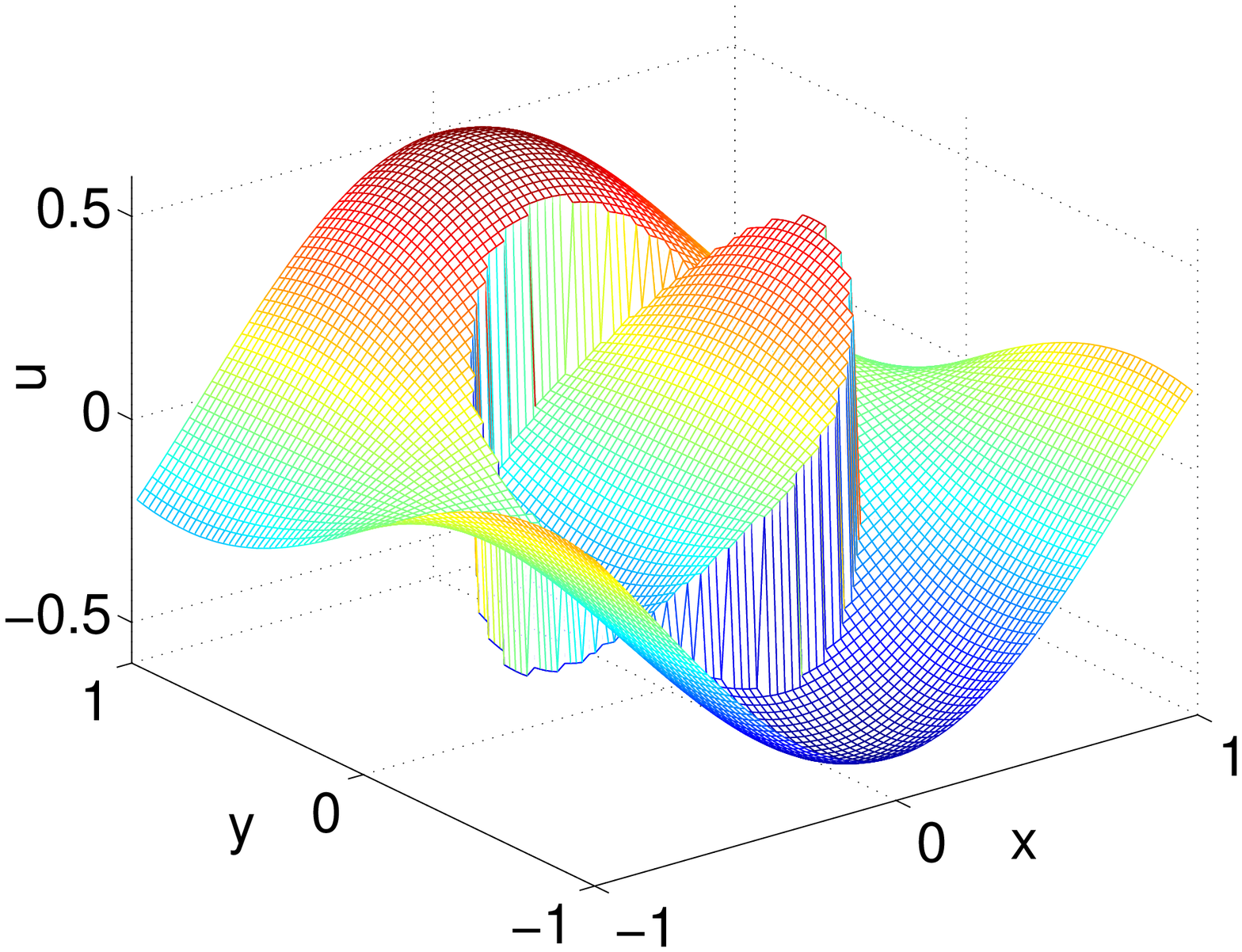,width=0.45\linewidth} \\
(a) &  (b) 
\end{tabular}
\end{center}
\caption{Numerical solution with $N=81$ at $T=1$. 
(a).  Example 3; (b). Example 4. 
}
\label{fig.solu34}
\end{figure*}

Again, the matched ADI algorithm is unconditionally stable for all
tested $h$ and $\Delta t$ in this example. By using $N=321$ and $T=1$, 
both $L_\infty$ and $L_2$ errors immediately begin to decay 
as $\Delta t$ becomes smaller,
and the limiting precisions are approached when $\Delta t$ is small enough. 
It can be seen from Fig. \ref{fig.ex34} (a) that the slope for the $L_\infty$ and $L_2$ 
error curve is, respectively, 1.07 and 1.06. Thus, the matched ADI delivers the first
order of accuracy in time for parabolic interface problems with general jump values. 
By taking $\Delta t=10^{-4}$, the spatial errors are reported in 
Table \ref{table.ex34}. The matched ADI method clearly achieves second order 
in space for this example. 
The matched ADI solution based on a mesh $N=81$ at the time $T=1$ is
shown in Fig. \ref{fig.solu34} (a). The jump values clearly change with respect to
the angle $\theta$.

\begin{table}[!t]
\caption{Spatial convergence tests for Example 3 and Example 4. }
\label{table.ex34}
\begin{center}
\begin{tabular}{|l|l|l|l|l|l|l|l|l|}
\hline
& \multicolumn{4}{c|}{Example 3} & \multicolumn{4}{c|}{Example 4} \\
\cline{2-5} \cline{6-9}
& \multicolumn{2}{c|}{$L_\infty$} & \multicolumn{2}{c|}{$L^2$} 
& \multicolumn{2}{c|}{$L_\infty$} & \multicolumn{2}{c|}{$L^2$} \\
\cline{2-3} \cline{4-5} \cline{6-7} \cline{8-9}
$N$ & error & order & error & order & error & order & error & order  \\
\hline
21   &  9.12E-3 &         &  1.61E-3 &         & 4.77E-3 &         & 8.54E-4 & \\
41   &  2.51E-3 & 1.86 &  3.76E-4 & 2.10 & 1.32E-3 & 1.85 & 1.96E-4 & 2.12 \\
81   &  4.93E-4 & 2.35 &  7.24E-5 & 2.38 & 2.56E-4 & 2.37 & 3.57E-5 & 2.46 \\
161 &  7.47E-5 & 2.72 &  1.39E-5 & 2.38 & 3.91E-5 & 2.71 & 6.66E-6 & 2.42 \\
321 &  1.36E-5 & 2.45 &  3.04E-6 & 2.19 & 7.22E-6 & 2.44 & 1.57E-6 & 2.09\\
\hline
\end{tabular}
\end{center}
\end{table}

{\bf Example 4.}
We next extend the Example 3 to a more general situation, by considering
time dependent jumps. 
The parameters for the interface and domain are fixed to be $r=0.5$, $D=0.99$, 
$\alpha^- = 1$, and $\alpha^+=10$.
The analytical solution is constructed as 
\begin{equation*}
u(x,y,t)= \begin{cases}
\sin(kx) \cos(ky) \cos(t),  \quad   \mbox{ if }   r <0.5 \\
\cos(kx) \sin(ky) \cos(t),  \quad   \mbox{ if } r \ge 0.5,
\end{cases}
\end{equation*}
where the wavenumber is chosen as $k=2$. 
The source term is given as
\begin{equation*}
f(x,y,t)= \begin{cases}
(2k^2 \alpha^- \cos(t) -\sin(t))  \sin(kx)\cos(ky),  \quad    \mbox{ if }   r <0.5 \\
(2k^2 \alpha^+ \cos(t) -\sin(t))  \cos(kx)\sin(ky),  \quad \mbox{ if } r \ge 0.5.
\end{cases}
\end{equation*}
Now, the jump conditions at a interface point 
$(x,y)=(\frac{1}{2} \cos \theta, \frac{1}{2} \sin \theta)$ 
depend on both space and time
\begin{align*}
[u] =& \cos (\frac{k}{2} \cos \theta) \sin (\frac{k}{2} \sin \theta) \cos(t)-
\sin (\frac{k}{2} \cos \theta) \cos (\frac{k}{2} \sin \theta) \cos(t), \\
[\alpha u_n] =&
k \cos(t) (\alpha^- \sin \theta - \alpha^+ \cos \theta)
\sin (\frac{k}{2} \cos \theta) \sin (\frac{k}{2} \sin \theta) \\
+& k \cos(t) (\alpha^+ \sin \theta - \alpha^- \cos \theta)
\cos (\frac{k}{2} \cos \theta) \cos (\frac{k}{2} \sin \theta), \\
[u_{\tau}] =&  
k \cos(t) (\cos \theta + \sin \theta)
(\cos (\frac{k}{2} \cos \theta) \cos (\frac{k}{2} \sin \theta)
+ \sin (\frac{k}{2} \cos \theta) \sin (\frac{k}{2} \sin \theta) ).
\end{align*}
Such jump conditions are the most general ones for  parabolic
interface problems.

\begin{figure*}[!tb]
\begin{center}
\begin{tabular}{cc}
\psfig{figure=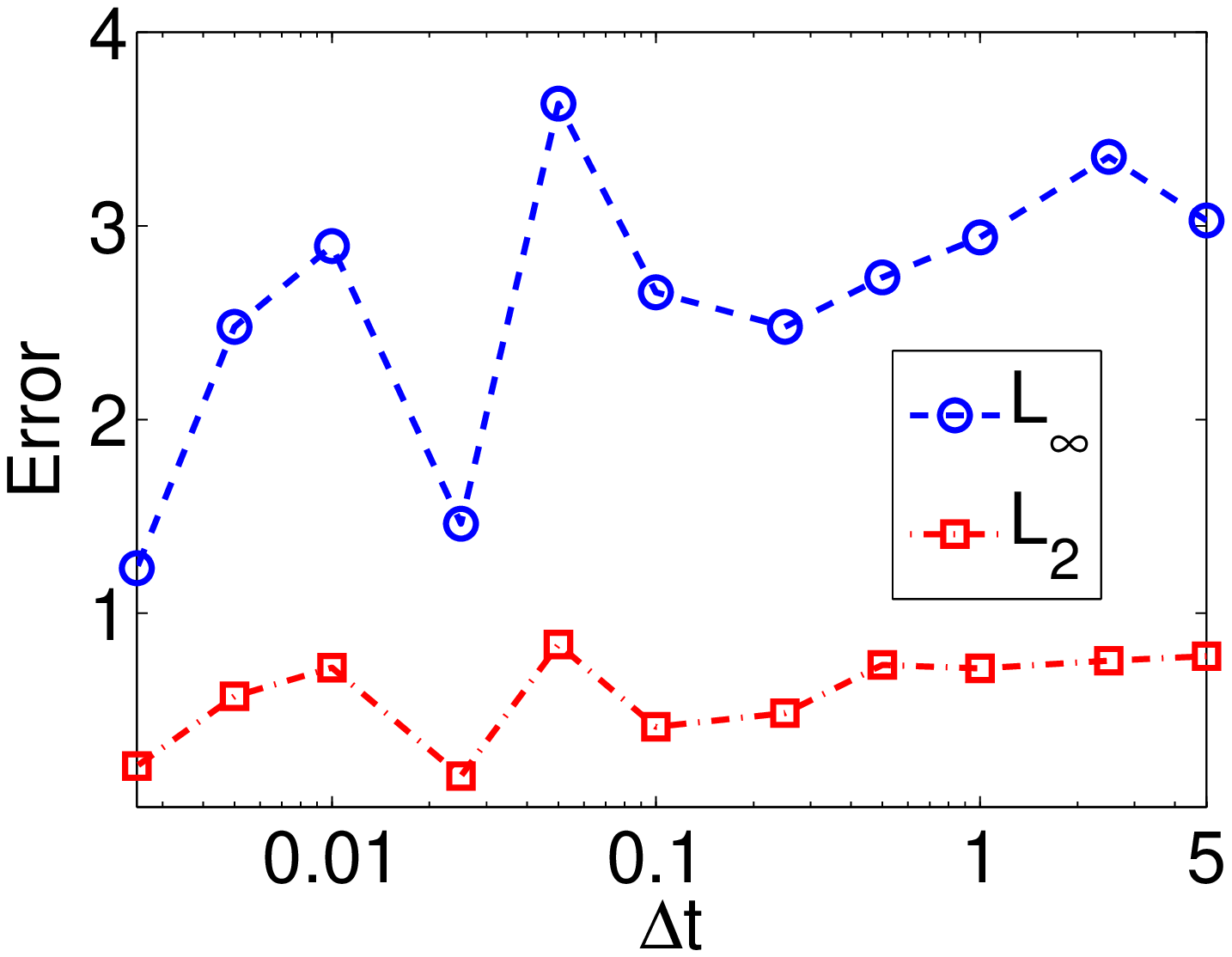,width=0.45\linewidth} & 
\psfig{figure=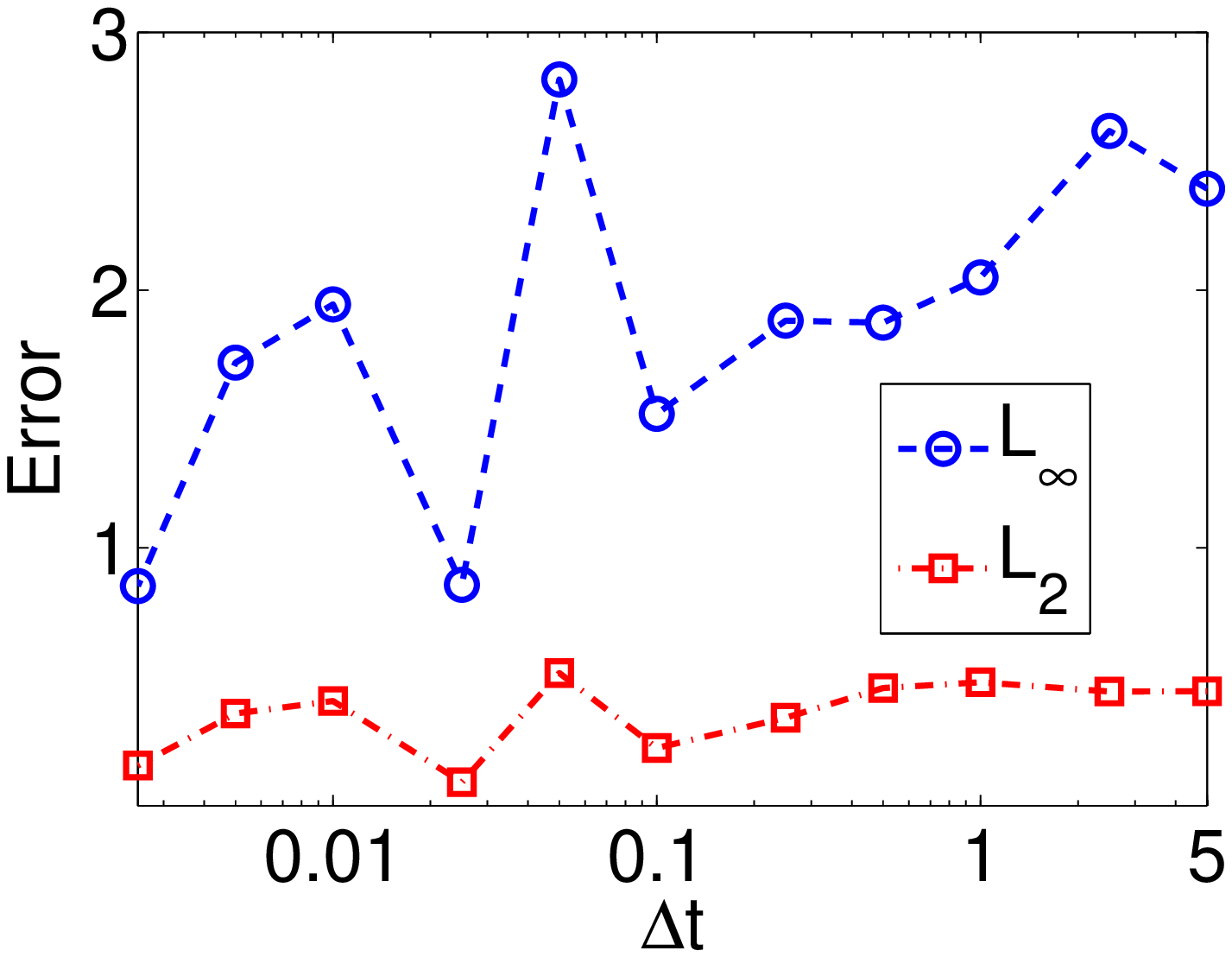,width=0.45\linewidth} \\
(a) &  (b) 
\end{tabular}
\end{center}
\caption{Bounded numerical errors with $N=321$
and $T=10^4 \Delta t$. 
(a).  Example 4; (b). Example 5, four leaves case.
}
\label{fig.large}
\end{figure*}

The temporal convergence pattern now turns out to be significantly different
from those of the previous examples. 
By using $N=321$ and $T=1$, the $L_\infty$ and $L_2$ errors
are depicted in Fig. \ref{fig.ex34} (b). 
It can be observed that the error curve does not immediately decay 
for large $\Delta t$ values. In other word, the temporal convergence
of the matched ADI method is somehow polluted by the time dependent
jump conditions. To relieve concerns about a potential instability for a
large $\Delta t$, the contaminated errors are depicted in Fig. \ref{fig.large} (a)
for $\Delta t$ values up to $\Delta t=5$.
For each $\Delta t$, we choose $N=321$ and $T=10^4 \Delta t$. 
This picture shows that
after $10^4$ time steps, all errors remain to be bounded. 
This demonstrates the unconditional stability of the matched ADI 
algorithm for solving time dependent jump conditions. 

As shown in Fig. \ref{fig.ex34} (b), only
when a rather small $\Delta t= 10^{-3}$ is employed,
the matched ADI method begins to converge. Nevertheless, once the 
convergence starts, the rate is pretty high. The least-squares fitting shows 
that the rate for the descending parts is $r=1.88$ and $r=1.91$, respectively,
for the $L_\infty$ and $L_2$ errors. 
Thus, excluding the polluted region, 
the matched ADI method yields a second order of accuracy in time. 
On the other hand, if the contaminated region was included in the least-squares
analysis, the overall temporal order would become about one.

The spatial orders are not affected by the time dependent jumps. 
By using a sufficiently small $\Delta t = 10^{-6}$, the numerical errors
based on different meshes are listed in Table \ref{table.ex34}. 
The numerical orders are all around two and are
comparable to those of the Example 3. 
This demonstrates the robustness of the proposed MIB interface 
treatment in solving time dependent jumps.
The matched ADI solution with $N=81$ and $T=1$ is plotted
in Fig \ref{fig.solu34} (b). The jump values shown in the figure will
oscillate with respect to the time $t$. 

\begin{figure*}[!tb]
\begin{center}
\begin{tabular}{cc}
\psfig{figure=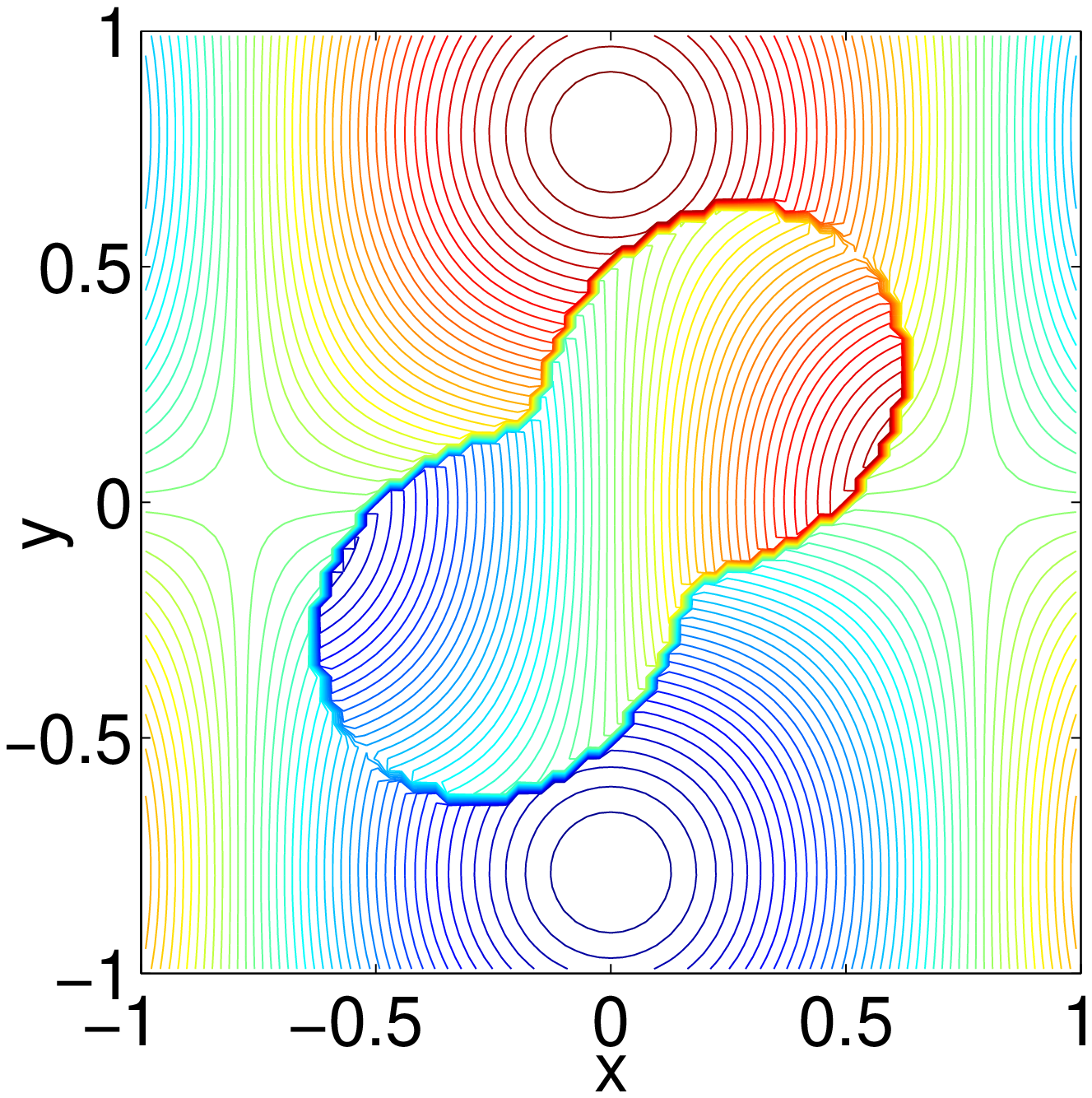,width=0.45\linewidth} & 
\psfig{figure=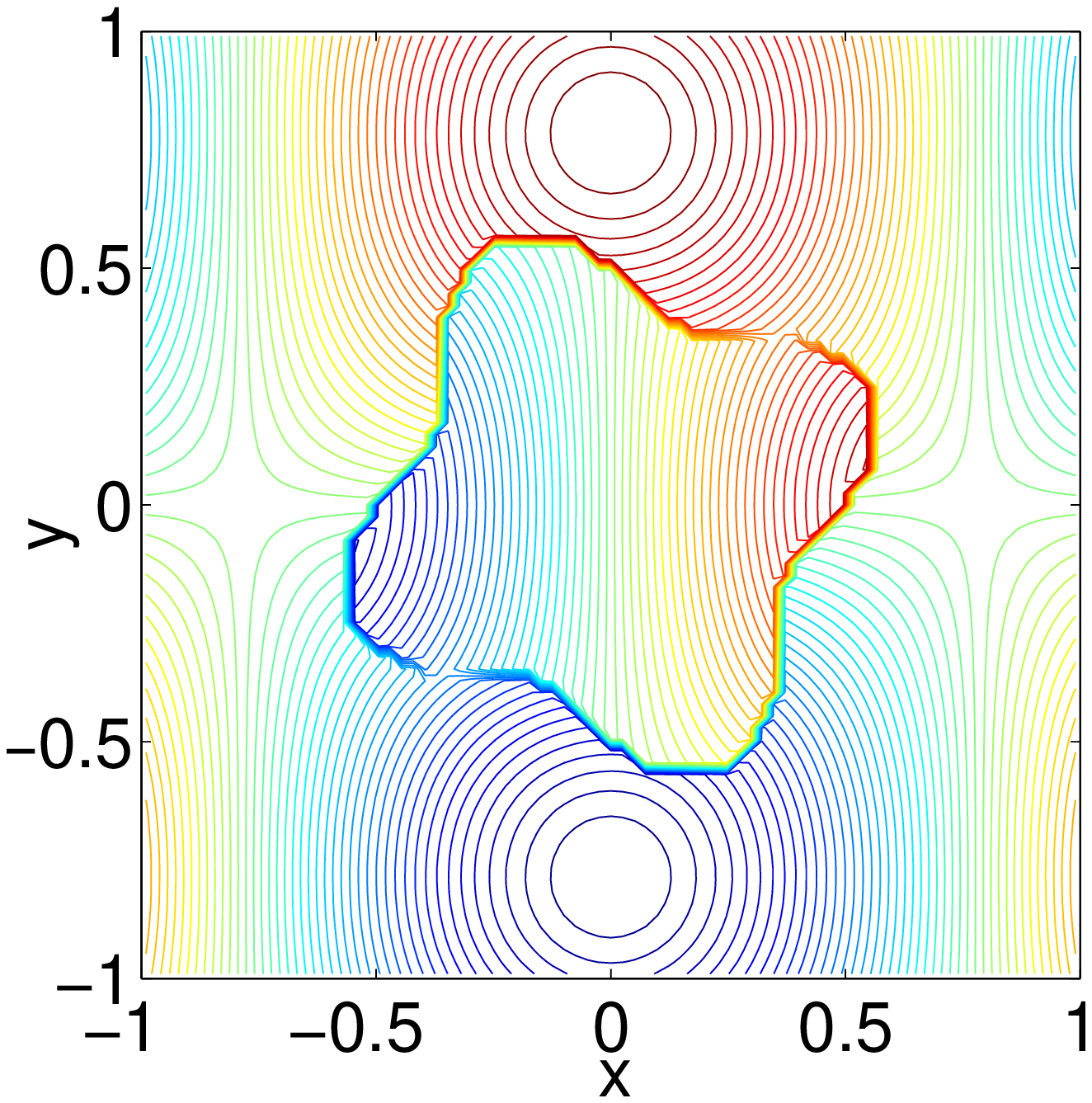,width=0.45\linewidth} \\
(a) &  (b) 
\end{tabular}
\end{center}
\caption{Contour plots of numerical solutions in Example 5. 
Here $N=81$ and $T=1$. 
(a).  Two leaves case;
(b). Four leaves case. 
}
\label{fig.solu5}
\end{figure*}

{\bf Example 5.}
At last, we explore the performance of the proposed matched ADI algorithm 
for interfaces of general shape. 
To this end, the following interface which is parameterized with the polar angle $s$
will be studied
\begin{equation}
\Gamma: \quad r= \frac{1}{2} + b \sin (m  s), \quad s \in [0,2 \pi]. 
\end{equation}
Here the parameter $m$ determines the number of ``leaves'' of the core region 
$\Omega^-$ and $b$ controls  the magnitude of the curvature.
Two independent cases with parameters $(m,b)=(2,1/4)$ and $(4,1/10)$
are considered. A square domain with $D=0.99$ is also employed. 
The resulting configurations of the two leaves and four leaves cases
can be seen from the contour plots of numerical 
solutions given in Fig. \ref{fig.solu5}. 
It is clear that concave segments or negative curvatures are involved in the 
present interfaces.
The analytical solution is constructed as in the Example 4
\begin{equation*}
u(x,y,t)= \begin{cases}
\sin(kx) \cos(ky) \cos(t),  \quad   \mbox{in } \Omega^-  \\
\cos(kx) \sin(ky) \cos(t),  \quad   \mbox{in } \Omega^+,
\end{cases}
\end{equation*}
with $k=2$. 
The source term is also given as
\begin{equation*}
f(x,y,t)= \begin{cases}
(2k^2 \alpha^- \cos(t) -\sin(t))  \sin(kx)\cos(ky),  \quad    \mbox{in } \Omega^- \\
(2k^2 \alpha^+ \cos(t) -\sin(t))  \cos(kx)\sin(ky),  \quad  \mbox{in } \Omega^+.
\end{cases}
\end{equation*}
The jump conditions can be similarly calculated according to the analytical solution.
The details are omitted here. 

We will rigorously examine the stability of the proposed matched ADI method by
using the four leaves case. 
For this purpose, we numerically calculate the leading eigenvalues of the
magnifying matrix ${\bf M}={\bf D}^{-1} {\bf B}$ with the largest magnitudes.
For a given $h$ and $\Delta t$, both ${\bf D}$ and ${\bf B}$ are saved in
a sparse matrix format. The inverse of  ${\bf D}$ is carried out by using a 
biconjugate gradient iterative solver, 
while the eigenvalues are computed by the eigenvalue package ARPACK. 
The tolerance is set to be $10^{-14}$ in these algebraic solvers. 
Because there are usually multiple leading eigenvalues whose magnitudes 
are the same, we will report the largest ten eigenvalues in magnitude, 
instead of just one spectral radius. Without the loss of generality, we denote
these ten eigenvalues to be $\lambda_i$ with $i=1,\ldots,10$ and
$|\lambda_i| \ge |\lambda_{i+1}|$. 
The spatial-temporal discretization can be
claimed to be stable, if the magnitudes of these leading eigenvalues are all
less than or equal to one, i.e., $|\lambda_i| \le 1$. 

\begin{figure*}[!tb]
\begin{center}
\begin{tabular}{cc}
\psfig{figure=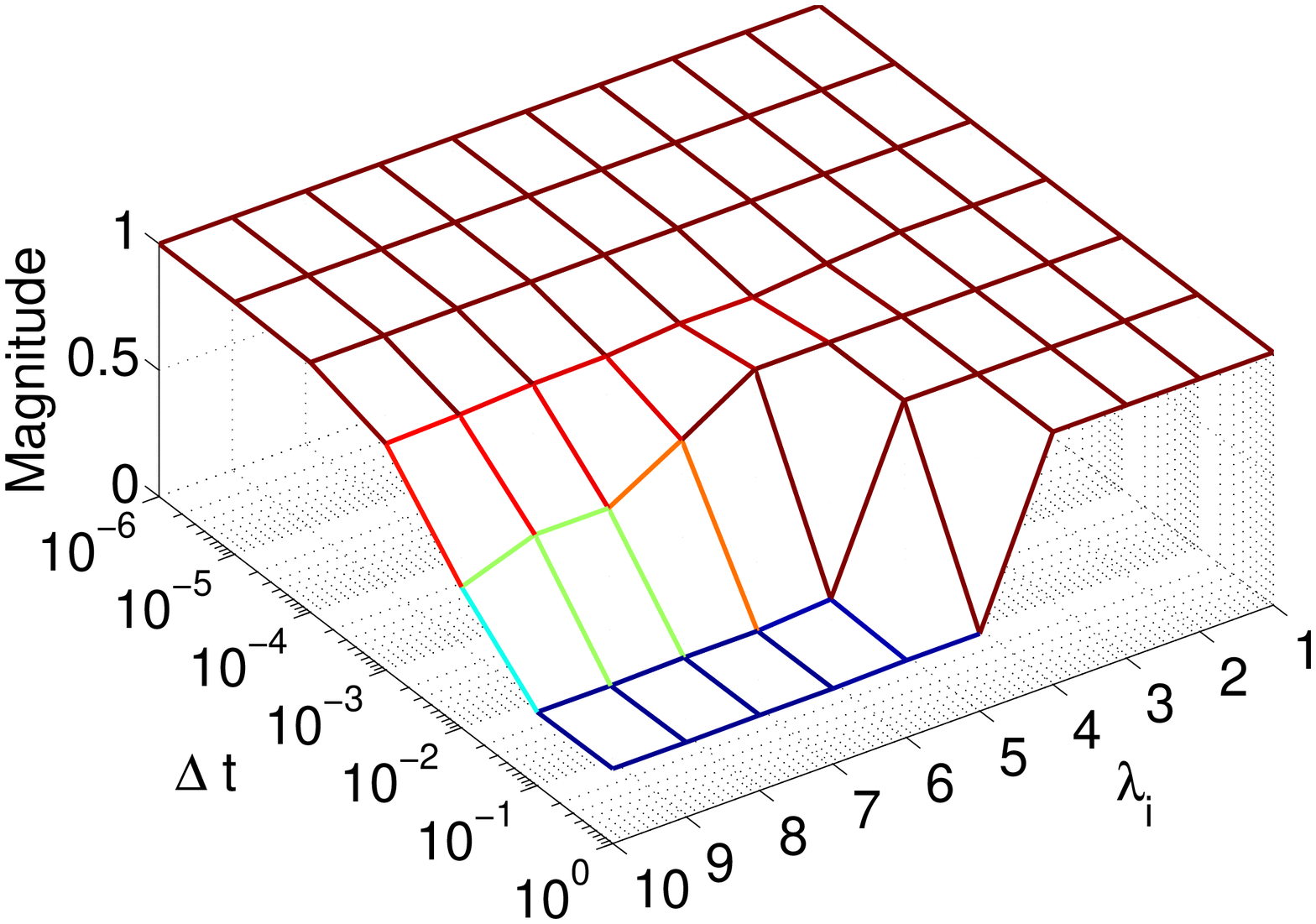,width=0.45\linewidth} & 
\psfig{figure=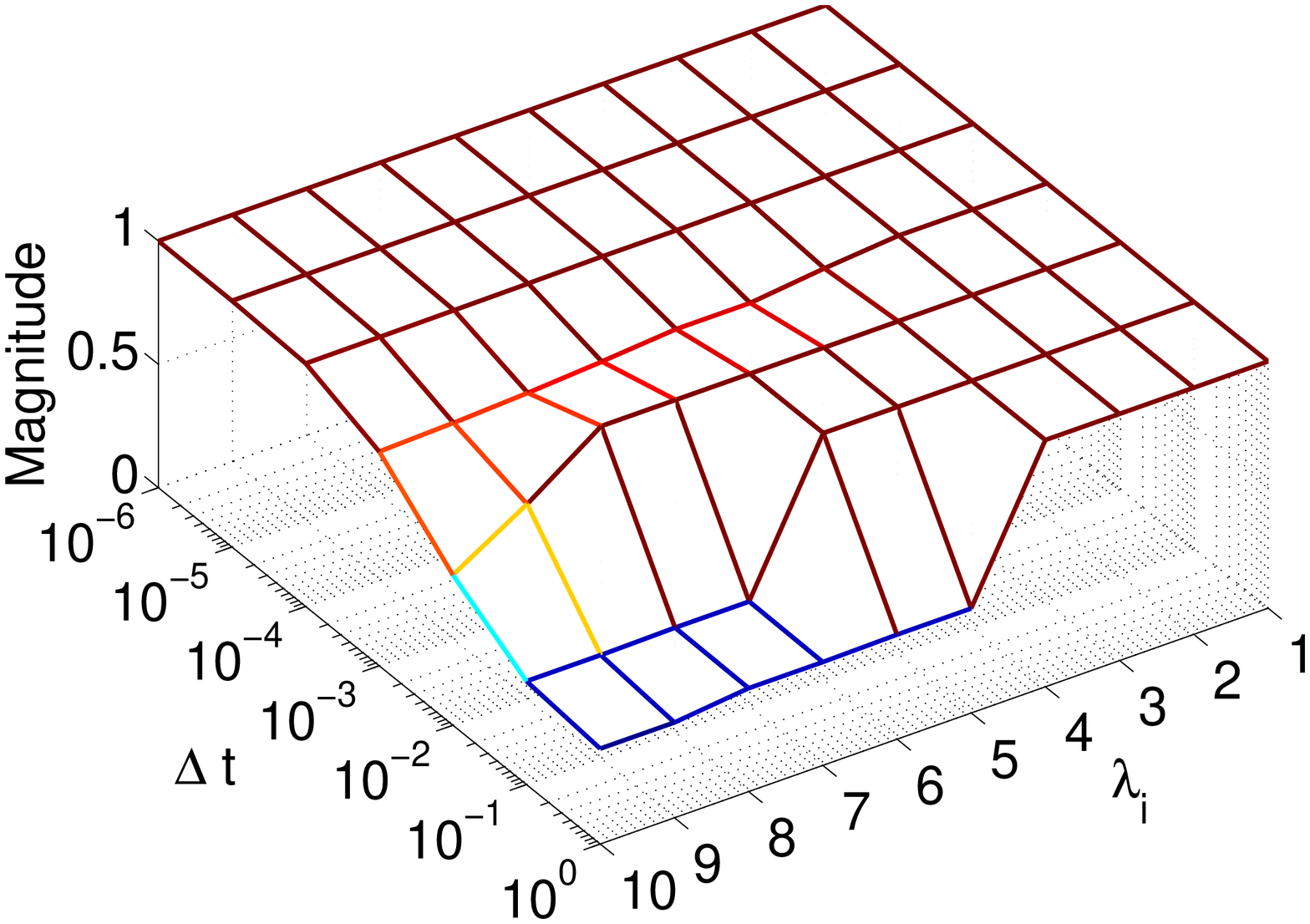,width=0.45\linewidth} \\
(a) &  (b) 
\end{tabular}
\end{center}
\caption{Leading eigenvalues of ${\bf M}$ with $N=41$.
(a).  $\alpha^- = 1$ and  $\alpha^+=10$;
(b). $\alpha^- = 1$ and  $\alpha^+=1000$. 
}
\label{fig.eig}
\end{figure*}

We first study the impact of $\Delta t$ on the stability.
By taking $\alpha^- = 1$,  $\alpha^+=10$, and $N=41$, the leading
eigenvalues are shown in Fig. \ref{fig.eig} (a) for seven $\Delta t$ values.
For $\Delta t=1$, $\Delta t=0.1$, and $\Delta t=0.01$, respectively, 
there are four, five, and six eigenvalues whose magnitudes equal to one. 
The rest leading eigenvalues take much smaller magnitudes. However,
for $\Delta t \le 10^{-3}$, all ten $\lambda_i$ values have almost the
same height in Fig. \ref{fig.eig} (a). Actually, there are just three
$\lambda_i$ with $|\lambda_i|=1$. The magnitudes of other eigenvalues
are strictly less than one, but are very close to one. 
Since $|\lambda_i| \le 1$ in all cases, 
the corresponding ADI computations are always stable.

We then carry out the similar tests by considering a much larger 
$\alpha^+=1000$. The other parameters are chosen to be the same. 
It can be seen from Fig. \ref{fig.eig} (b) that the leading eigenvalues 
are very similar to those in Fig. \ref{fig.eig} (a). 
The only minor change is that 
for $\Delta t=1$, $\Delta t=0.1$, and $\Delta t=0.01$, respectively, 
there are four, six, and eight eigenvalues whose magnitudes equal to one. 
The present  analysis validates the stability of the 
matched ADI scheme for a large $\alpha^+$. 
This study also demonstrates the robustness of the MIB interface method in 
handling large jump ratios.

\begin{figure*}[!tb]
\begin{center}
\psfig{figure=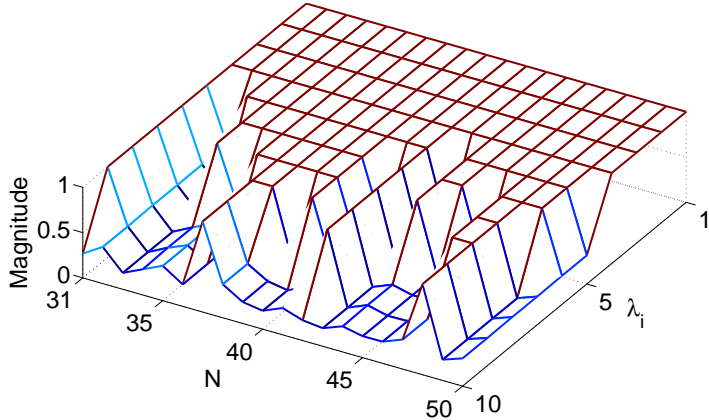,width=0.6\linewidth} 
\end{center}
\caption{Leading eigenvalues of ${\bf M}$ with $\Delta t=1$,
$\alpha^- = 1$, and  $\alpha^+=10$.
}
\label{fig.eig2}
\end{figure*}

We next calculate leading eigenvalues for a different $N$ with a fixed 
$\Delta t=1$. The results with $\alpha^- = 1$ and  $\alpha^+=10$ will be
reported, while those of $\alpha^+=1000$ are found to be similar. 
By considering 20 mesh sizes starting from $N=31$, the eigenvalues are 
shown in Fig. \ref{fig.eig2}. With a different $N$, the number of the
largest eigenvalues with the same magnitude
is at least three, and could sometimes be more than ten. 
In general, the dependence of this number with respect to $N$
is quite random. 
This agrees with our previous discussion that 
the geometry of a complicated interface $\Gamma$ will 
affect the spectrum of ${\bf M}$ in a random way, 
because the intersection of $\Gamma$ with a mesh of size $N$ by $N$ is 
quite arbitrary. With $|\lambda_i| \le 1$ for all $N$ values, 
the present study demonstrates that such a geometry effect will not 
compromise the unconditional stability of the matched ADI method.

\begin{figure*}[!tb]
\begin{center}
\begin{tabular}{cc}
\psfig{figure=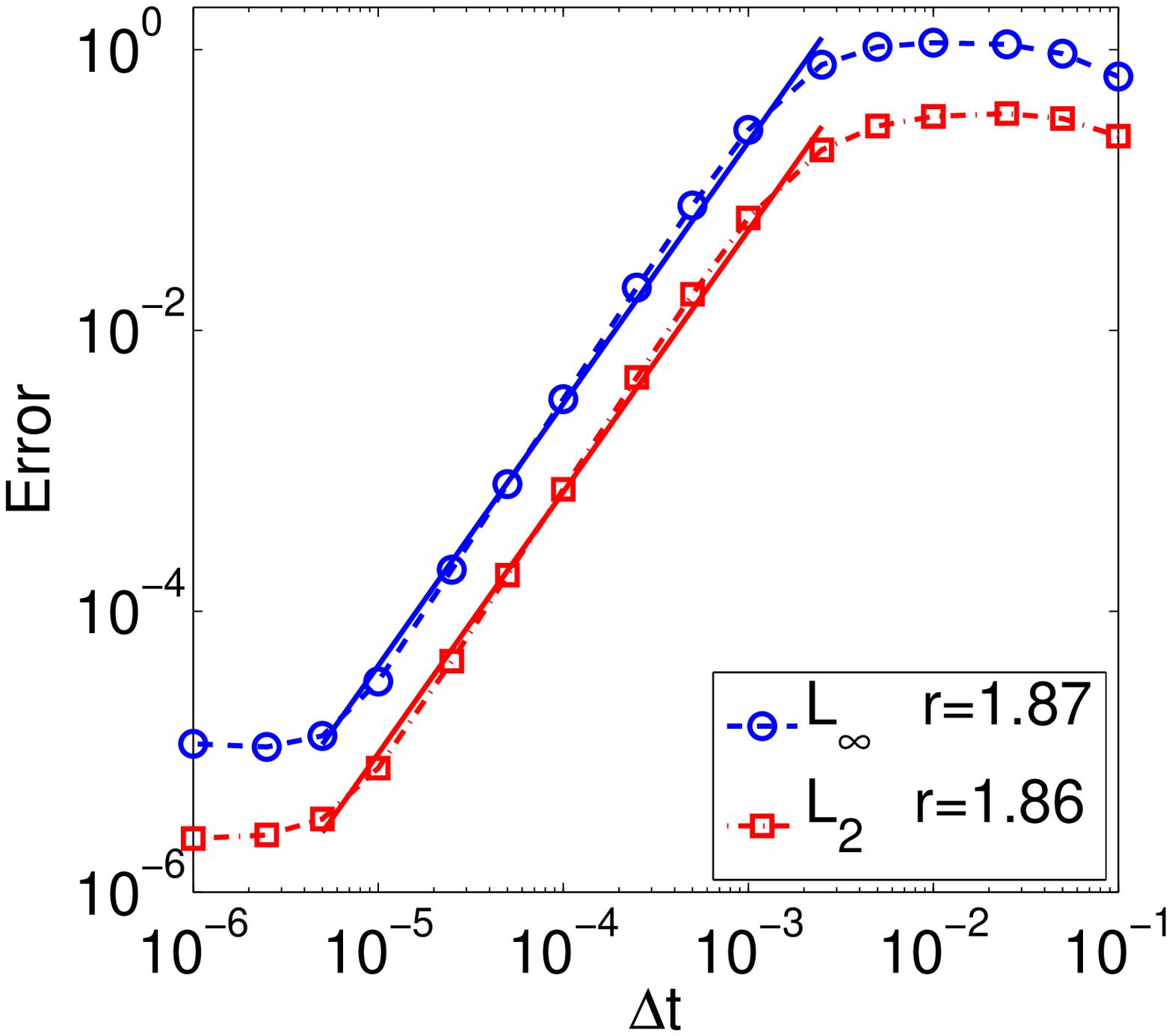,width=0.45\linewidth} & 
\psfig{figure=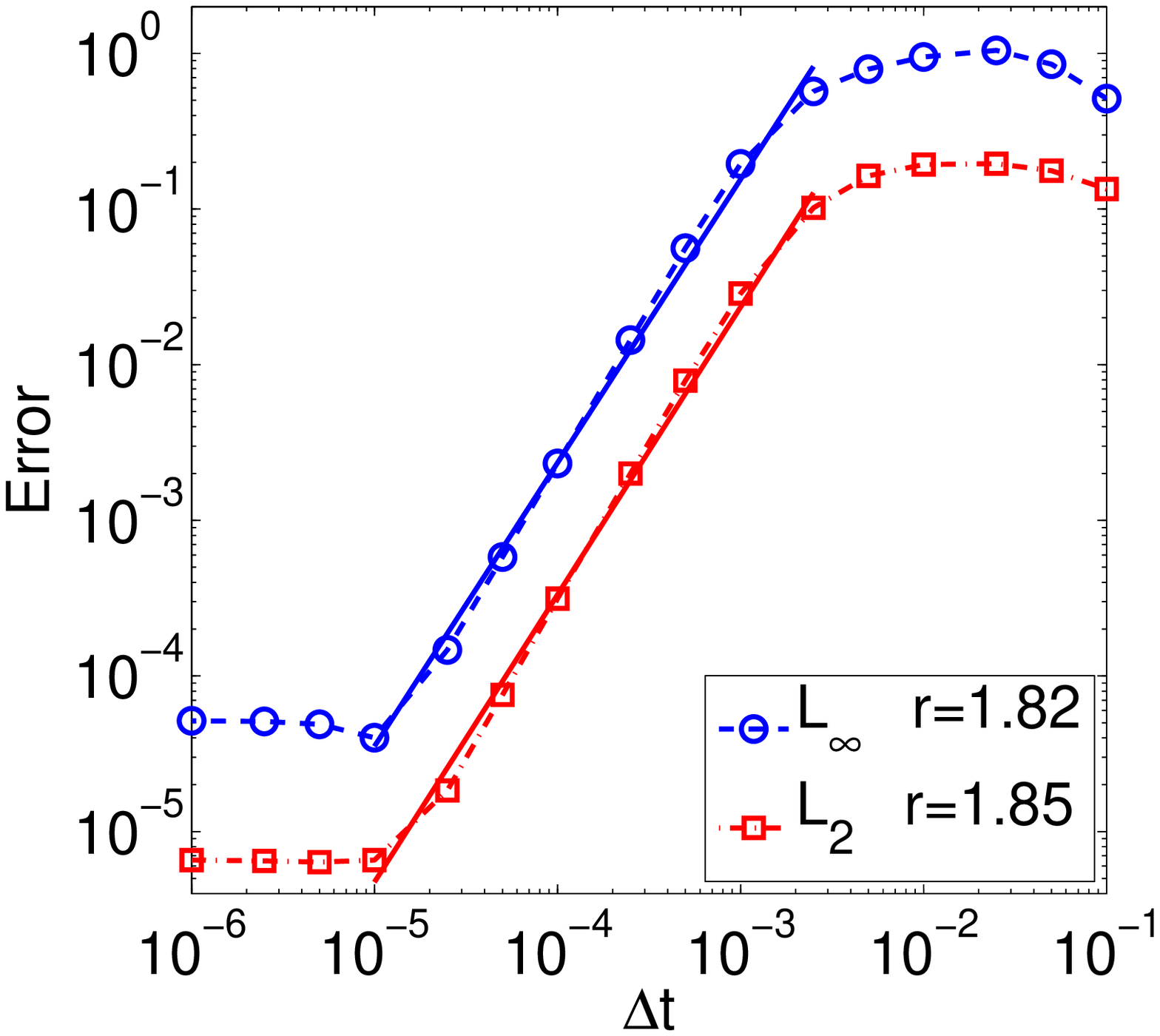,width=0.45\linewidth} \\
(a) &  (b) 
\end{tabular}
\end{center}
\caption{Temporal convergence tests of Example 5. 
(a).  Two leaves case;
(b). Four leaves case. 
}
\label{fig.ex5}
\end{figure*}

After establishing the stability of the matched ADI scheme, 
we next examine the temporal  convergence by considering 
$\alpha^-=1$ and $\alpha^+=10$. 
By using $N=321$, $T=1$ and different $\Delta t$ values,
the $L_\infty$ and $L_2$ errors for both
two leaves and four leaves cases are shown in Fig. \ref{fig.ex5}. 
The temporal convergences of these two cases are very similar to that of Example 4,
because the present jump conditions are also time dependent. 
In particular, the convergence is polluted for large $\Delta t$ values. 
Similarly, the contaminated errors remain to be bounded for a 
long time-stepping, see Fig \ref{fig.large}. (b). 
As shown in Fig. \ref{fig.ex5},
the temporal convergence begins only when $\Delta t$ is small enough 
and is of second order once occurs. 
The overall temporal order would be also around one, if the polluted region was also 
included in the least-squares analysis.

\begin{table}[!t]
\caption{Spatial convergence tests for Example 5. }
\label{table.ex5}
\begin{center}
\begin{tabular}{|l|l|l|l|l|l|l|l|l|}
\hline
& \multicolumn{4}{c|}{Two leaves case} & \multicolumn{4}{c|}{Four leaves case} \\
\cline{2-5} \cline{6-9}
& \multicolumn{2}{c|}{$L_\infty$} & \multicolumn{2}{c|}{$L^2$} 
& \multicolumn{2}{c|}{$L_\infty$} & \multicolumn{2}{c|}{$L^2$} \\
\cline{2-3} \cline{4-5} \cline{6-7} \cline{8-9}
$N$ & error & order & error & order & error & order & error & order  \\
\hline
21   &  3.06E-3 &         &  7.49E-4 &         & 5.97E-3 &         & 1.48E-3 & \\
41   &  5.37E-4 & 2.51 &  1.47E-4 & 2.35 & 2.48E-3 & 1.27 & 5.83E-4 & 1.34 \\
81   &  1.80E-4 & 1.58 &  4.59E-5 & 1.68 & 9.95E-4 & 1.32 & 1.63E-4 & 1.84 \\
161 &  3.92E-5 & 2.20 &  1.03E-5 & 2.16 & 1.32E-4 & 2.91 & 1.86E-5 & 3.13 \\
321 &  1.08E-5 & 1.86 &  2.55E-6 & 2.01 & 4.01E-5 & 1.72 & 6.57E-6 & 1.50\\
\hline
\end{tabular}
\end{center}
\end{table}

We finally study the spatial convergence with $\alpha^-=1$ and $\alpha^+=10$. 
By taking $\Delta t=2.5 \times 10^{-6}$ and $\Delta t= 10^{-5}$, respectively, for the
two leaves case and four leaves case, the spatial errors are reported in 
Table \ref{table.ex5}. 
It can be observed  that when the geometrical structure
becomes more complicated or the number of leaves $m$ is larger, 
the convergence pattern of the matched ADI method becomes more oscillatory 
and the overall numerical order becomes slightly smaller. 
Nevertheless, the proposed matched ADI method still can secure a second order
of accuracy for these challenging parabolic interface problems of complicated
geometry. 
The contour plots of the matched ADI solutions based on a mesh $N=81$ 
at time $T=1$ are illustrated in Fig. \ref{fig.solu5}. The solutions clearly undergo
sharp changes across the interface $\Gamma$, and such changes are time
variant.

\section{Conclusion}\label{Sec:Conclusion}

This paper presents a novel matched ADI method for solving parabolic interface 
problems with general jump conditions and complex geometries. 
The development of accurate and stable ADI schemes for such interface problems
is essentially an open problem, because the existing IIM-ADI schemes 
require the second order jump conditions and have difficulties to correct finite difference
approximations based on general flux jumps. 
The second order jump conditions are not needed in our interface treatment.
Moreover, a novel tensor product decomposition is proposed
to decouple 2D jump conditions into essentially 1D ones. By enforcing these 1D conditions,
rigorous 1D MIB schemes are developed to treat regular interfaces and corner interfaces.
The resulting matched ADI scheme achieves second order of accuracy in space 
and first order of accuracy in time for interfaces of different shapes. 
The efficiency of the ADI scheme is well maintained, because the MIB interface treatment
needs to be conducted only once at the beginning of the computation and fast algebraic
solvers are developed for perturbed tridiagonal systems.
Stability analysis by means of the numerical spectrum analysis of the magnifying
matrix is conducted to examine the impact of geometry and various parameters.
The matched ADI scheme is found to be unconditional stable with all 
numerical eigenvalues having magnitude less than or equal to one. 
However, the stability proof of the matched ADI method
remains to be an open question, 
because the finite difference weights of the MIB discretization 
depend on the interface geometry in an unpredictable manner.

In our preliminary studies, 
the direct application of the present MIB spatial discretization with the 
Peaceman-Rachford ADI method is found to be conditionally stable. 
The development of robust matched ADI schemes with second order in time and
for more general parabolic equations is currently under our investigation.

\vspace{1cm}

\centerline{\bf Acknowledgment}
\noindent
This work was supported in part by NSF grants DMS-1016579 and DMS-1318898,
and the University of Alabama Research Stimulation Program (RSP) award.


\begin{thebibliography}{99}

\bibitem{AdamsLi} L.M. Adams and Z.L. Li,
The immersed interface/multigrid method for interface problems,
\emph{SIAM J. Sci. Comput.}, {\bf 24}, 463-479, (2002).

\bibitem{Attanayake}  C. Attanayake and D. Senaratne,
Convergence of an immersed finite element method for semilinear
parabolic interface problems,
\emph{Appl. Math. Sci.}, {\bf 5}, 135-147, (2011).

\bibitem{Bates09} P.W. Bates, Z. Chen, Y.H. Sun, G.W. Wei, and S. Zhao,
Geometric and potential driving formation and evolution of biomolecular surfaces,
\emph{J. Math. Biol.}, {\bf 59}, 193-231, (2009).

\bibitem{Bouchon10a} F. Bouchon and G.H. Peichl,
An immersed interface technique for the numerical solution of the heat equation on 
a moving domain, 
In \emph{Numerical Mathematics and Advanced Applications 2009}, 
G. Kreiss et al. (eds.), Springer-Verlag Berlin Heidelberg, 
181-189, (2010). 

\bibitem{Bouchon10b}  F. Bouchon and G.H. Peichl,
The immersed interface technique for parabolic problems with mixed
boundary conditions,
\emph{SIAM J. Numer. Anal.}, {\bf 48}, 2247-2266, (2010).


\bibitem{ChenZou} Z. Chen and J. Zou,
Finite element methods and their convergence for elliptic and parabolic 
interface problems,
\emph{Numer. Math.}, {\bf 79}, 175-202, (1996).

\bibitem{Chen11} L. Chen, J. Shen, and C. Xu,
Spectral direction spitting schemes for the incompressible Navier-Stokes
equations, 
\emph{East Asia J. Appl. Math.}, {\bf 1}, 215-234, (2011).


\bibitem{Douglas} J. Douglas, Jr.,
On the numerical integration of 
$\frac{\partial^2 u}{\partial x^2} + \frac{\partial^2 u}{\partial y^2} 
= \frac{\partial u}{\partial t}$ by implicit methods,  
\emph{J. Soc. Indust. Appl. Math.}, {\bf 3}, 42-65, (1955).

\bibitem{DougPeace} J. Douglas, Jr. and D. Peaceman,
Numerical solution of two-dimensional heat flow problems, 
\emph{Amer. Inst. Chem. Engrg. J.}, {\bf 1}, 505-512, (1955).

\bibitem{Geng13} W. Geng and S. Zhao,
Fully implicit ADI schemes for solving the nonlinear Poisson-Boltzmann equation,
\emph{Molecular Based Mathematical Biology}, {\bf 1}, 109-123, (2013).

\bibitem {Kandilarov04}  J.D. Kandilarov and L.G. Vulkov,
The immersed interface method for a nonlinear chemical diffusion
equation with local sites of reactions, 
\emph{Numer. Algorithms}, {\bf 36}, 2850307, (2004).

\bibitem {Kandilarov07}  J.D. Kandilarov and L.G. Vulkov,
The immersed interface method for two-dimensional heat-diffusion
equations with singular own sources, 
\emph{Appl. Numer. Math.}, {\bf 57}, 486-497, (2007).


\bibitem{LevLi} R.J. LeVeque and Z.L. Li,
The immersed interface method for elliptic equations with
discontinuous coefficients and singular sources,
\emph{SIAM J. Numer. Anal.}, {\bf 31}, 1019-1044, (1994).

\bibitem{LiMayo} Z.L. Li and A. Mayo,
ADI methods for heat equations with discontinuities along an arbitrary
interface, In \emph{Proc. Symp. Appl. Math.}, (W. Gautschi ed.) AMS, 
{\bf 48}, 311-315, (1993).

\bibitem{LiShen} Z.L. Li and Y.-Q. Shen,
A numerical method for solving heat equations involving interfaces, 
\emph{Electron. J. Diff. Eqns.}, Conf. 03, 100-108, (2000).

\bibitem{Liu13}  J. Liu and Z. Zheng,
IIM-based ADI finite difference scheme for nonlinear 
convection-diffusion equations with interfaces, 
\emph{Appl. Math. Modelling}, {\bf 37}, 1196-1207, (2013).

\bibitem{Liu14}  J. Liu and Z. Zheng,
A dimension by dimension splitting immersed interface method
for heat conduction equation with interfaces, 
\emph{J. Comput. Appl. Math.}, {\bf 261}, 221-231, (2014).


\bibitem{PeaceRach}  D. Peaceman and H. Rachford,
The numerical solution of parabolic and elliptic equations, 
\emph{J. Soc. Indust. Appl. Math.}, {\bf 3}, 28-41, (1955).

\bibitem{NR}  W.H. Press, S.A. Teukolsky, W.T. Vetterling, and B.P. Flannery,
\emph{Numerical Recipes: The Art of Scientific Computing} (3rd ed.)
New York: Cambridge University Press, (2007). 

\bibitem{Sinha05}  R.K. Sinha and B. Deka,
Optimal error estimates for linear parabolic problems with discontinuous
coefficients, 
\emph{SIAM J. Numer. Anal.}, {\bf 43}, 733-749, (2005).

\bibitem{Sinha09}  R.K. Sinha and B. Deka,
Finite element methods for semi linear elliptic and parabolic interface
problems, 
\emph{Appl. Numer. Math.}, {\bf 59}, 1870-1883, (2009).


\bibitem{ADIbook} J.C. Strikwerda, \emph{Finite difference schemes
and partial differential equations}, SIAM, 2nd edition, (2004).

\bibitem {Tay} W.C. Tay, D.Y. Heh, and E.L. Tan,
GPU-accelerated fundamental ADI-FDTD with complex frequency shifted
convolutional perfectly matched layer, 
\emph{Progress In Electromagnetics Research}, {\bf 14}, 177-192, (2010).


\bibitem {Tian14} W. Tian and S. Zhao,
A fast ADI algorithm for geometric flow equations in biomolecular surface generation,
\emph{Int. J. Numer. Meth. Biomed. Engrg.}, {\bf 30}, 490-516, (2014).

\bibitem {Wang10} S. Wang, R. Samulyak, and T. Guo,
An embedded boundary method for elliptic and parabolic problems with 
interfaces and application to multi-material systems with phase transitions,
\emph{Acta Math. Sci. Ser. B}, {\bf 30}, 499-521, (2010).

\bibitem {Wei13} Z. Wei, B. Jang, Y. Zhang, and Y. Jia,
Parallelizing alternating direction implicit solver on GPUs, 
\emph{Procedia Computer Science}, {\bf 18}, 389-398, (2013).

\bibitem {Zhao07} S. Zhao,
On the spurious solutions in the high-order finite difference 
methods, 
\emph{Comput. Method Appl. Mech. Engrg.}, {\bf 196}, 5031-5046, (2007).

\bibitem {Zhao14} S. Zhao,
Operator splitting ADI schemes for
pseudo-time coupled nonlinear solvation simulations,
\emph{J. Comput. Phys.}, {\bf 257}, 1000-1021, (2014).

\bibitem {Zhao04} S. Zhao and G.W. Wei,
High order FDTD methods via derivative matching for Maxwell's equations
with material interfaces, 
\emph{J. Comput. Phys.}, {\bf 200}, 60-103, (2004).

\bibitem{Zhao09} S. Zhao and G.W. Wei,
Matched interface and boundary (MIB) method for the implementation of
boundary conditions in high-order central finite differences,
\emph{Int. J. Numer. Method Engrg.}, {\bf 77}, 1690-1730, (2009).


\bibitem{Zhou06} Y.C. Zhou, S. Zhao, M. Feig, and G.W. Wei,
High order matched interface and boundary method for elliptic equations
with discontinuous coefficients and singular sources, 
\emph{J. Comput. Phys.}, {\bf 213}, 1-30, (2006).

\end{thebibliography}
\end{document}